%
\newif\ifloadreferences\loadreferencestrue
%
%
%
%
%
\let\myfrac=\frac%
\input eplain %
\let\frac=\myfrac%
\let\myfootnote=\footnote%
\input amstex \input epsf %
\let\footnote=\myfootnote%
%
%
\loadeufm\loadmsam\loadmsbm\message{symbol names}\UseAMSsymbols\message{,}%
\magnification 1200 %
\font\myfontdefault=cmr10%
\newif\ifmakebiblio%
\newif\ifinappendices%
\newif\ifundefinedreferences%
\newif\ifchangedreferences%
\makebibliofalse%
\undefinedreferencesfalse%
\changedreferencesfalse%
%
%
%
%
%
\def\setcatcodes{\catcode`\!=0 \catcode`\\=11}%
{\global\let\noe=\noexpand%
\catcode`\@=11 \catcode`\_=11 \setcatcodes%
!global!def!_@@internal@@makeref#1{%
!global!expandafter!def!csname #1ref!endcsname##1{%
!csname _@#1@##1!endcsname%
!expandafter!ifx!csname _@#1@##1!endcsname!relax%
    !write16{#1 ##1 not defined - run saving references}%
    !undefinedreferencestrue%
!fi}}%
!global!def!_@@internal@@makelabel#1{%
!global!expandafter!def!csname #1label!endcsname##1{%
!edef!temptoken{!csname #1info!endcsname}%
!ifloadreferences%
    !expandafter!ifx!csname _@#1@##1!endcsname!relax%
        !write16{#1 ##1 not hitherto defined - rerun saving references}%
        !changedreferencestrue%
    !else%
        !expandafter!ifx!csname _@#1@##1!endcsname!temptoken%
        !else%
            !write16{#1 ##1 reference has changed - rerun saving references}%
            !changedreferencestrue%
        !fi%
    !fi%
!else%
    !expandafter!edef!csname _@#1@##1!endcsname{!temptoken}%
    !edef!textoutput{!write!references{\global\def\_@#1@##1{!temptoken}}}%
    !textoutput%
!fi}}%
!global!def!makecounter#1{!_@@internal@@makelabel{#1}!_@@internal@@makeref{#1}}%
!unsetcatcodes%
}
%
%
%
%
%
\def\turnintolatin#1{\ifcase #1 _\or i\or ii\or iii\or iv\or v\or vi\or vii\or viii\or ix\or x\or xi\or xii\or xiii\or xiv\or xv\or xvi\or xvii\or xviii\or xix\or xx\or xxi\or xxii\or xxiii\or xxiv\or xxv\or xxvi\fi}%
\def\alphanum#1{\ifcase #1 _\or A\or B\or C\or D\or E\or F\or G\or H\or I\or J\or K\or L\or M\or N\or O\or P\or Q\or R\or S\or T\or U\or V\or W\or X\or Y\or Z\fi}%
\newwrite\references%
\ifloadreferences{\catcode`\@=11 \catcode`\_=11 \global\def\_@citation@BrendleEichmair{1}
\global\def\_@citation@Chavel{2}
\global\def\_@citation@EichmairMetzgerI{3}
\global\def\_@citation@EichmairMetzgerII{4}
\global\def\_@citation@EichmairMetzgerIII{5}
\global\def\_@citation@Krylov{6}
\global\def\_@citation@MaxNunSmi{7}
\global\def\_@citation@NardulliI{8}
\global\def\_@citation@NardulliII{9}
\global\def\_@citation@Salamon{10}
\global\def\_@citation@RosSmi{11}
\global\def\_@citation@Rudin{12}
\global\def\_@citation@Schwarz{13}
\global\def\_@citation@SmiCCH{14}
\global\def\_@citation@SmiEMCFI{15}
\global\def\_@citation@White{16}
\global\def\_@citation@Ye{17}
\global\def\_@head@Introduction{1}
\global\def\_@subhead@Background{1.1}
\global\def\_@subhead@NotationTerminologyAndMainResult{1.2}
\global\def\_@eqn@GradientFlowEquation{\relax \unhbox \voidb@x \hbox {{\relax \tenrm (1)}}}
\global\def\_@eqn@LinearisedGradientFlowEquation{\relax \unhbox \voidb@x \hbox {{\relax \tenrm (2)}}}
\global\def\_@eqn@AreaMinusVolumeFunctional{\relax \unhbox \voidb@x \hbox {{\relax \tenrm (3)}}}
\global\def\_@eqn@ForcedMeanCurvatureFlowEquation{\relax \unhbox \voidb@x \hbox {{\relax \tenrm (4)}}}
\global\def\_@eqn@PerturbedFlowOfSpheres{\relax \unhbox \voidb@x \hbox {{\relax \tenrm (5)}}}
\global\def\_@proc@ThmMainTheorem{1.2.1}
\global\def\_@subhead@Discussion{1.3}
\global\def\_@eqn@ScaleDependenceOperator{\relax \unhbox \voidb@x \hbox {{\relax \tenrm (6)}}}
\global\def\_@subhead@OverviewAcknowledgements{1.4}
\global\def\_@head@TheTaylorSeriesOfGeometricFunctions{2}
\global\def\_@subhead@CurvatureTensors{2.1}
\global\def\_@eqn@MetricIsExponential{\relax \unhbox \voidb@x \hbox {{\relax \tenrm (7)}}}
\global\def\_@eqn@SpacesOfCurvatureTensors{\relax \unhbox \voidb@x \hbox {{\relax \tenrm (8)}}}
\global\def\_@proc@TheSpaceOfCurvaturePolynomialsIsSelfAdjoint{2.1.1}
\global\def\_@eqn@DefnOfOverlineR{\relax \unhbox \voidb@x \hbox {{\relax \tenrm (9)}}}
\global\def\_@subhead@CurvaturePolynomials{2.2}
\global\def\_@eqn@FundamentalCurvaturePolynomial{\relax \unhbox \voidb@x \hbox {{\relax \tenrm (10)}}}
\global\def\_@eqn@SpaceOfCurvaturePolynomials{\relax \unhbox \voidb@x \hbox {{\relax \tenrm (11)}}}
\global\def\_@proc@CurvaturePolynomialsInOneVariableHaveNonNegativeOrderDifference{2.2.1}
\global\def\_@eqn@ExamplesOfSpacesOfCurvaturePolynomials{\relax \unhbox \voidb@x \hbox {{\relax \tenrm (12)}}}
\global\def\_@eqn@MoreExamplesOfSpacesOfCurvaturePolynomials{\relax \unhbox \voidb@x \hbox {{\relax \tenrm (13)}}}
\global\def\_@eqn@ExtendedCurvaturePolynomials{\relax \unhbox \voidb@x \hbox {{\relax \tenrm (14)}}}
\global\def\_@eqn@ExamplesOfSpacesOfExtendedCurvaturePolynomials{\relax \unhbox \voidb@x \hbox {{\relax \tenrm (15)}}}
\global\def\_@subhead@GeneralPropertiesOfTaylorSeries{2.3}
\global\def\_@eqn@BinomialTheorem{\relax \unhbox \voidb@x \hbox {{\relax \tenrm (16)}}}
\global\def\_@eqn@DefinitionOfPowerOfSeriesByBinomialTheorem{\relax \unhbox \voidb@x \hbox {{\relax \tenrm (17)}}}
\global\def\_@proc@PowersOfCurvatureSeriesAreAlsoCurvatureSeries{2.3.1}
\global\def\_@proc@TestForBeingCurvatureSeries{2.3.2}
\global\def\_@subhead@TensorValuedGeometricFunction{2.4}
\global\def\_@eqn@EqnParallelTransport{\relax \unhbox \voidb@x \hbox {{\relax \tenrm (18)}}}
\global\def\_@eqn@TaylorSeriesOfM{\relax \unhbox \voidb@x \hbox {{\relax \tenrm (19)}}}
\global\def\_@proc@TaylorSeriesOfM{2.4.1}
\global\def\_@eqn@EqnMetric{\relax \unhbox \voidb@x \hbox {{\relax \tenrm (20)}}}
\global\def\_@eqn@EqnShortAsymptoticSeriesForAAndB{\relax \unhbox \voidb@x \hbox {{\relax \tenrm (21)}}}
\global\def\_@eqn@TaylorSeriesOfAAndB{\relax \unhbox \voidb@x \hbox {{\relax \tenrm (22)}}}
\global\def\_@proc@TaylorSeriesOfAAndB{2.4.2}
\global\def\_@eqn@EqnChristoffelSymbol{\relax \unhbox \voidb@x \hbox {{\relax \tenrm (23)}}}
\global\def\_@eqn@EqnShortSeriesForChristoffelSymbol{\relax \unhbox \voidb@x \hbox {{\relax \tenrm (24)}}}
\global\def\_@eqn@TaylorSeriesOfChristoffelSymbol{\relax \unhbox \voidb@x \hbox {{\relax \tenrm (25)}}}
\global\def\_@proc@TaylorSeriesOfChristoffelSymbol{2.4.3}
\global\def\_@eqn@KoszulFormula{\relax \unhbox \voidb@x \hbox {{\relax \tenrm (26)}}}
\global\def\_@subhead@TheExponentialMapAndParallelTransport{2.5}
\global\def\_@eqn@DomainOfExponentialMap{\relax \unhbox \voidb@x \hbox {{\relax \tenrm (27)}}}
\global\def\_@eqn@TaylorSeriesOfExp{\relax \unhbox \voidb@x \hbox {{\relax \tenrm (28)}}}
\global\def\_@eqn@FirstTermsOfTaylorSeriesOfExp{\relax \unhbox \voidb@x \hbox {{\relax \tenrm (29)}}}
\global\def\_@proc@TaylorSeriesOfExp{2.5.1}
\global\def\_@eqn@TaylorSeriesOfT{\relax \unhbox \voidb@x \hbox {{\relax \tenrm (30)}}}
\global\def\_@eqn@FirstTermsOfTaylorSeriesOfT{\relax \unhbox \voidb@x \hbox {{\relax \tenrm (31)}}}
\global\def\_@proc@TaylorSeriesOfParallelTransport{2.5.2}
\global\def\_@eqn@AdvancedDomainOfExponentialMap{\relax \unhbox \voidb@x \hbox {{\relax \tenrm (32)}}}
\global\def\_@eqn@InductiveDefinitionOfExpAndTStepI{\relax \unhbox \voidb@x \hbox {{\relax \tenrm (33)}}}
\global\def\_@eqn@InductiveDefinitionOfExpAndTStepII{\relax \unhbox \voidb@x \hbox {{\relax \tenrm (34)}}}
\global\def\_@eqn@TaylorSeriesOfExtendedFunctions{\relax \unhbox \voidb@x \hbox {{\relax \tenrm (35)}}}
\global\def\_@eqn@FirstTermsOfTaylorSeriesOfExtendedFunctions{\relax \unhbox \voidb@x \hbox {{\relax \tenrm (36)}}}
\global\def\_@proc@TaylorSeriesOfExtendedFunctions{2.5.3}
\global\def\_@head@TaylorSeriesOfFunctionsDerivedFromImmersions{3}
\global\def\_@subhead@GraphsOverSpheres{3.1}
\global\def\_@eqn@GraphOverSphere{\relax \unhbox \voidb@x \hbox {{\relax \tenrm (37)}}}
\global\def\_@eqn@DefinitionsOfXAndR{\relax \unhbox \voidb@x \hbox {{\relax \tenrm (38)}}}
\global\def\_@eqn@FormulaForF{\relax \unhbox \voidb@x \hbox {{\relax \tenrm (39)}}}
\global\def\_@eqn@GradientOfF{\relax \unhbox \voidb@x \hbox {{\relax \tenrm (40)}}}
\global\def\_@proc@GradientOfF{3.1.1}
\global\def\_@eqn@FormulaeForRAndY{\relax \unhbox \voidb@x \hbox {{\relax \tenrm (41)}}}
\global\def\_@eqn@FormulaForNHat{\relax \unhbox \voidb@x \hbox {{\relax \tenrm (42)}}}
\global\def\_@eqn@FormulaForN{\relax \unhbox \voidb@x \hbox {{\relax \tenrm (43)}}}
\global\def\_@eqn@FormulaForImmersionVariableCentre{\relax \unhbox \voidb@x \hbox {{\relax \tenrm (44)}}}
\global\def\_@subhead@TheTaylorSeriesOfTheUnitNormalVector{3.2}
\global\def\_@eqn@GeneratorsOfCurvaturePolynomials{\relax \unhbox \voidb@x \hbox {{\relax \tenrm (45)}}}
\global\def\_@eqn@ExtendedCurvaturePolynomialsOfTheSecondKind{\relax \unhbox \voidb@x \hbox {{\relax \tenrm (46)}}}
\global\def\_@eqn@MeaningOfFormalism{\relax \unhbox \voidb@x \hbox {{\relax \tenrm (47)}}}
\global\def\_@eqn@TaylorSeriesOfROverT{\relax \unhbox \voidb@x \hbox {{\relax \tenrm (48)}}}
\global\def\_@proc@TaylorSeriesOfROverT{3.2.1}
\global\def\_@eqn@TaylorSeriesOfNormN{\relax \unhbox \voidb@x \hbox {{\relax \tenrm (49)}}}
\global\def\_@eqn@ShortTaylorSeriesOfNormN{\relax \unhbox \voidb@x \hbox {{\relax \tenrm (50)}}}
\global\def\_@proc@TaylorSeriesOfNormN{3.2.2}
\global\def\_@eqn@TaylorSeriesOfNPartI{\relax \unhbox \voidb@x \hbox {{\relax \tenrm (51)}}}
\global\def\_@eqn@TaylorSeriesOfNPartII{\relax \unhbox \voidb@x \hbox {{\relax \tenrm (52)}}}
\global\def\_@eqn@FirstTermsInTaylorSeriesOfN{\relax \unhbox \voidb@x \hbox {{\relax \tenrm (53)}}}
\global\def\_@proc@TaylorSeriesOfN{3.2.3}
\global\def\_@eqn@TaylorSeriesOfNVariableCentrePartI{\relax \unhbox \voidb@x \hbox {{\relax \tenrm (54)}}}
\global\def\_@eqn@TaylorSeriesOfNVariableCentrePartII{\relax \unhbox \voidb@x \hbox {{\relax \tenrm (55)}}}
\global\def\_@eqn@FirstTermsInTaylorSeriesOfNVariableCentre{\relax \unhbox \voidb@x \hbox {{\relax \tenrm (56)}}}
\global\def\_@proc@TaylorSeriesOfNVariableCentre{3.2.4}
\global\def\_@subhead@NormalVariation{3.3}
\global\def\_@eqn@MostGeneralVariation{\relax \unhbox \voidb@x \hbox {{\relax \tenrm (57)}}}
\global\def\_@eqn@DefnOfPQR{\relax \unhbox \voidb@x \hbox {{\relax \tenrm (58)}}}
\global\def\_@eqn@DefnOfpqr{\relax \unhbox \voidb@x \hbox {{\relax \tenrm (59)}}}
\global\def\_@eqn@TaylorSeriesOfLittlePPartI{\relax \unhbox \voidb@x \hbox {{\relax \tenrm (60)}}}
\global\def\_@eqn@TaylorSeriesOfLittleP{\relax \unhbox \voidb@x \hbox {{\relax \tenrm (61)}}}
\global\def\_@eqn@FirstTermsInTaylorSeriesOfLittleP{\relax \unhbox \voidb@x \hbox {{\relax \tenrm (62)}}}
\global\def\_@proc@TaylorSeriesOfLittleP{3.3.1}
\global\def\_@eqn@TaylorSeriesOfLittleQPartI{\relax \unhbox \voidb@x \hbox {{\relax \tenrm (63)}}}
\global\def\_@eqn@TaylorSeriesOfLittleQ{\relax \unhbox \voidb@x \hbox {{\relax \tenrm (64)}}}
\global\def\_@eqn@FirstTermsInTaylorSeriesOfLittleQ{\relax \unhbox \voidb@x \hbox {{\relax \tenrm (65)}}}
\global\def\_@proc@TaylorSeriesOfLittleQ{3.3.2}
\global\def\_@eqn@TaylorSeriesOfLittleR{\relax \unhbox \voidb@x \hbox {{\relax \tenrm (66)}}}
\global\def\_@eqn@FirstTermsInTaylorSeriesOfLittleR{\relax \unhbox \voidb@x \hbox {{\relax \tenrm (67)}}}
\global\def\_@proc@TaylorSeriesOfLittleR{3.3.3}
\global\def\_@subhead@TheTaylorSeriesOfTheMeanCurvature{3.4}
\global\def\_@eqn@FirstAsymptoticFormulaForDisplacedCurvature{\relax \unhbox \voidb@x \hbox {{\relax \tenrm (68)}}}
\global\def\_@eqn@TaylorSeriesOfHPartI{\relax \unhbox \voidb@x \hbox {{\relax \tenrm (69)}}}
\global\def\_@eqn@TaylorSeriesOfHPartII{\relax \unhbox \voidb@x \hbox {{\relax \tenrm (70)}}}
\global\def\_@proc@TaylorSeriesOfH{3.4.1}
\global\def\_@eqn@FormulaForMeanCurvature{\relax \unhbox \voidb@x \hbox {{\relax \tenrm (71)}}}
\global\def\_@eqn@GradientOfFAgain{\relax \unhbox \voidb@x \hbox {{\relax \tenrm (72)}}}
\global\def\_@eqn@SecondPartOfMeanCurvature{\relax \unhbox \voidb@x \hbox {{\relax \tenrm (73)}}}
\global\def\_@head@AsymptoticExpansionsAndFormalSolutions{4}
\global\def\_@subhead@AsymptoticExpansions{4.1}
\global\def\_@eqn@DefinitionOfAsymptoticExpansionPartI{\relax \unhbox \voidb@x \hbox {{\relax \tenrm (74)}}}
\global\def\_@eqn@DefinitionOfAsymptoticExpansionPartII{\relax \unhbox \voidb@x \hbox {{\relax \tenrm (75)}}}
\global\def\_@eqn@AsymptoticSeriesForP{\relax \unhbox \voidb@x \hbox {{\relax \tenrm (76)}}}
\global\def\_@proc@AsymptoticSeriesForP{4.1.1}
\global\def\_@eqn@AsymptoticSeriesForQ{\relax \unhbox \voidb@x \hbox {{\relax \tenrm (77)}}}
\global\def\_@proc@AsymptoticSeriesForQ{4.1.2}
\global\def\_@eqn@AsymptoticSeriesForR{\relax \unhbox \voidb@x \hbox {{\relax \tenrm (78)}}}
\global\def\_@proc@AsymptoticSeriesForR{4.1.3}
\global\def\_@eqn@AsymptoticExpansionOfH{\relax \unhbox \voidb@x \hbox {{\relax \tenrm (79)}}}
\global\def\_@eqn@FirstTermsInAsymptoticExpansionOfH{\relax \unhbox \voidb@x \hbox {{\relax \tenrm (80)}}}
\global\def\_@proc@AsymptoticSeriesForH{4.1.4}
\global\def\_@subhead@FlowsOfSurfaces{4.2}
\global\def\_@eqn@FamilyOfImmersions{\relax \unhbox \voidb@x \hbox {{\relax \tenrm (81)}}}
\global\def\_@eqn@VariationOfImmersions{\relax \unhbox \voidb@x \hbox {{\relax \tenrm (82)}}}
\global\def\_@eqn@ForcedMCFEquation{\relax \unhbox \voidb@x \hbox {{\relax \tenrm (83)}}}
\global\def\_@eqn@AsymptoticExpansionOfPhi{\relax \unhbox \voidb@x \hbox {{\relax \tenrm (84)}}}
\global\def\_@eqn@FirstTermsInAsymptoticSeriesOfPhi{\relax \unhbox \voidb@x \hbox {{\relax \tenrm (85)}}}
\global\def\_@proc@AsymptoticExpansionOfPhi{4.2.1}
\global\def\_@subhead@ParabolicOperatorsI{4.3}
\global\def\_@eqn@DefinitionOfHolderSeminorm{\relax \unhbox \voidb@x \hbox {{\relax \tenrm (86)}}}
\global\def\_@eqn@DefinitionOfHolderNorm{\relax \unhbox \voidb@x \hbox {{\relax \tenrm (87)}}}
\global\def\_@eqn@DefinitionOfHolderSpace{\relax \unhbox \voidb@x \hbox {{\relax \tenrm (88)}}}
\global\def\_@eqn@DefinitionOfOperatorP{\relax \unhbox \voidb@x \hbox {{\relax \tenrm (89)}}}
\global\def\_@eqn@DefinitionOfPEpsilon{\relax \unhbox \voidb@x \hbox {{\relax \tenrm (90)}}}
\global\def\_@eqn@DefinitionOfGEpsilon{\relax \unhbox \voidb@x \hbox {{\relax \tenrm (91)}}}
\global\def\_@eqn@FirstWeightedHolderNorm{\relax \unhbox \voidb@x \hbox {{\relax \tenrm (92)}}}
\global\def\_@eqn@OperatorNormOfGreensOperator{\relax \unhbox \voidb@x \hbox {{\relax \tenrm (93)}}}
\global\def\_@proc@NormOfInverseOfFiniteDimensionalParabolicOperator{4.3.1}
\global\def\_@subhead@ParabolicOperatorsII{4.4}
\global\def\_@eqn@DefinitionOfInhomogeneousHolderSeminorms{\relax \unhbox \voidb@x \hbox {{\relax \tenrm (94)}}}
\global\def\_@eqn@DefinitionOfInhomogeneousHolderNorm{\relax \unhbox \voidb@x \hbox {{\relax \tenrm (95)}}}
\global\def\_@eqn@DefinitionOfInhomogeneousHolderSpace{\relax \unhbox \voidb@x \hbox {{\relax \tenrm (96)}}}
\global\def\_@eqn@DefinitionOfWeightedHolderNorm{\relax \unhbox \voidb@x \hbox {{\relax \tenrm (97)}}}
\global\def\_@eqn@DefinitionOfOperatorQ{\relax \unhbox \voidb@x \hbox {{\relax \tenrm (98)}}}
\global\def\_@eqn@FunctionsInOrthogonalComplement{\relax \unhbox \voidb@x \hbox {{\relax \tenrm (99)}}}
\global\def\_@eqn@InfiniteDimensionalFactor{\relax \unhbox \voidb@x \hbox {{\relax \tenrm (100)}}}
\global\def\_@proc@InfiniteDimensionalFactor{4.4.1}
\global\def\_@eqn@NonLinearSphericalHarmonics{\relax \unhbox \voidb@x \hbox {{\relax \tenrm (101)}}}
\global\def\_@proc@NormOfRestrictionOfGreensOperator{4.4.2}
\global\def\_@subhead@MoreOnSphericalHarmonics{4.5}
\global\def\_@eqn@IsotropicTensors{\relax \unhbox \voidb@x \hbox {{\relax \tenrm (102)}}}
\global\def\_@eqn@SymmetricProduct{\relax \unhbox \voidb@x \hbox {{\relax \tenrm (103)}}}
\global\def\_@proc@SymmetricIsotropicTensors{4.5.1}
\global\def\_@eqn@Contraction{\relax \unhbox \voidb@x \hbox {{\relax \tenrm (104)}}}
\global\def\_@eqn@ContractionOfProductWithDelta{\relax \unhbox \voidb@x \hbox {{\relax \tenrm (105)}}}
\global\def\_@proc@ContractionOfProductWithDelta{4.5.2}
\global\def\_@eqn@ContractionOfPowerWithDelta{\relax \unhbox \voidb@x \hbox {{\relax \tenrm (106)}}}
\global\def\_@proc@ContractionOfPowerWithDelta{4.5.3}
\global\def\_@eqn@IntegralFormulae{\relax \unhbox \voidb@x \hbox {{\relax \tenrm (107)}}}
\global\def\_@proc@IntegralFormulae{4.5.4}
\global\def\_@proc@OrthonormalBasis{4.5.5}
\global\def\_@eqn@ZerothTermInPhi{\relax \unhbox \voidb@x \hbox {{\relax \tenrm (108)}}}
\global\def\_@eqn@ThirdSummand{\relax \unhbox \voidb@x \hbox {{\relax \tenrm (109)}}}
\global\def\_@proc@ThirdSummand{4.5.6}
\global\def\_@eqn@FirstTermInPhi{\relax \unhbox \voidb@x \hbox {{\relax \tenrm (110)}}}
\global\def\_@eqn@FourthSummand{\relax \unhbox \voidb@x \hbox {{\relax \tenrm (111)}}}
\global\def\_@proc@FourthSummand{4.5.7}
\global\def\_@subhead@FormalSolutions{4.6}
\global\def\_@eqn@SequenceProperties{\relax \unhbox \voidb@x \hbox {{\relax \tenrm (112)}}}
\global\def\_@eqn@SequencesYieldFormalSolution{\relax \unhbox \voidb@x \hbox {{\relax \tenrm (113)}}}
\global\def\_@proc@ExistenceOfApproximateSolutions{4.6.1}
\global\def\_@eqn@QMapsToKernelOfPi{\relax \unhbox \voidb@x \hbox {{\relax \tenrm (114)}}}
\global\def\_@subhead@ExactSolutions{4.7}
\global\def\_@proc@ThmExistence{4.7.2}
\global\def\_@eqn@SingularDerivative{\relax \unhbox \voidb@x \hbox {{\relax \tenrm (115)}}}
\global\def\_@head@Bibliography{5}
 }%
\else{\openout\references=references.tex }%
\fi%
%
%
\newcount\headno%
\global\headno=0%
\def\headinfo{\ifinappendices\alphanum\headno\else\the\headno\fi}%
\def\nextheadno[#1]{\global\advance\headno by 1 \global\subheadno=0 \global\procno=0 \headinfo\headlabel{#1}}%
\makecounter{head}%
%
%
\newcount\subheadno%
\global\subheadno=0%
\def\subheadinfo{\headinfo.\the\subheadno}%
\def\nextsubheadno[#1]{\global\advance\subheadno by 1 \global\procno=0 \subheadinfo\subheadlabel{#1}}%
\makecounter{subhead}%
%
%
\newcount\procno%
\global\procno=0%
\def\procinfo{\subheadinfo.\the\procno}%
\def\nextprocno{\global\advance\procno by 1 \procinfo}%
\makecounter{proc}%
%
%
\newcount\figno%
\global\figno=0%
\def\figinfo{\subheadinfo.\the\figno}%
\def\nextfigno{\global\advance\figno by 1 \figinfo}%
\makecounter{fig}%
%
%
\newcount\eqnno%
\global\eqnno=0%
\def\eqninfo{\text{{\rm (\the\eqnno)}}}%
\def\nexteqnno[#1]{\global\advance\eqnno by 1\hbox to 0em{\eqnlabel{#1}}\eqninfo}%
\makecounter{eqn}%
%
%
%
%
%
\def\gobbleeight#1#2#3#4#5#6#7#8{}%
\newcount\citationno%
\global\citationno=0%
\def\citationinfo{\the\citationno}%
\makecounter{citation}%
\newwrite\biblio%
\def\newref#1#2{%
\def\temptext{#2}%
\edef\bibliotextoutput{\expandafter\gobbleeight\meaning\temptext}%
\global\advance\citationno by 1\citationlabel{#1}%
\ifmakebiblio%
    \edef\fileoutput{\write\biblio{\noindent\hbox to 0pt{\hss$[\the\citationno]$}\hskip 0.2em\bibliotextoutput\medskip}}%
    \fileoutput%
\fi}%
\def\cite#1{%
$[\citationref{#1}]$%
\ifmakebiblio%
    \edef\fileoutput{\write\biblio{#1}}%
    \fileoutput%
\fi%
}%
%
%
%
%
\let\mypar=\par%
\edef\Pagetitle={Blank}\headline={\hfil\Pagetitle\hfil}%
\edef\Pagefooter={Blank}\footline={\hfil\Pagefooter\hfil}%
%
%
\newcount\showpagenumflag%
\global\showpagenumflag=0 %
\def\nextoddpage%
{\newpage\ifodd\pageno%
\else\global\showpagenumflag=0 %
\null\vfil\eject%
\global\showpagenumflag=1 %
\fi}%
%
%
\font\headfont=cmb12%
\def\newhead#1[#2]%
{\ifhmode\mypar\fi%
\ifnum\headno=0 \else\goodbreak\bigskip\fi%
{\headfont\noindent\nextheadno[#2]\ - #1.}
\nobreak\medskip}%
%
%
\def\newsubhead#1[#2]%
{\ifhmode\mypar\fi%
\ifnum\subheadno=0 \else\goodbreak\medskip\fi%
{\bf\noindent\nextsubheadno[#2]\ - #1.\ }}%
%
%
\newif\ifinproclaim%
\global\inproclaimfalse%
\def\proclaim#1{%
\goodbreak\medskip
\bgroup\inproclaimtrue%
\noindent{\bf #1}%
\nobreak\medskip\sl}%
\def\noskipproclaim#1{%
\goodbreak\medskip%
\bgroup\inproclaimtrue%
\noindent{\bf #1}\nobreak\sl}%
\def\endproclaim{\mypar\egroup\nobreak\medskip\ignorespaces}%
%
%
%
\newcount\xpos\newcount\ypos
\def\makelabelgrid{%
\xpos=-5 \ypos=-5 %
\loop\ifnum\xpos<6 %
{\loop\ifnum\ypos<6 %
\def\labeltext{x}%
\ifnum\xpos=0\def\labeltext{+}\fi%
\ifnum\ypos=0\def\labeltext{+}\fi%
\placelabel[\xpos][\ypos]{\labeltext}%
\advance\ypos by 1 %
\repeat}%
\advance\xpos by 1 %
\repeat}%
\def\placelabel[#1][#2]#3{{%
\setbox10=\hbox{\raise #2cm \hbox{\hskip #1cm #3}}%
\ht10=0pt \dp10=0pt \wd10=0pt \box10}}%
%
%
%
%
\def\myitem#1{\noindent\hbox to .5cm{\hfill#1\hss}}%
%
%
%
%
%
%
%
%
%
\font\sansseriften=cmss10%
\font\sansserifseven=cmss7%
\font\sansseriffive=cmss5%
\newfam\sansseriffam%
\textfont\sansseriffam=\sansseriften%
\scriptfont\sansseriffam=\sansserifseven%
\scriptscriptfont\sansseriffam=\sansseriffive%
\def\mathsf{\fam\sansseriffam}%
%
%
%
\font\boldten=cmb10%
\font\boldseven=cmb7%
\font\boldfive=cmb5%
\newfam\mathboldfam%
\textfont\mathboldfam=\boldten%
\scriptfont\mathboldfam=\boldseven%
\scriptscriptfont\mathboldfam=\boldfive%
\def\mathbf{\fam\mathboldfam}%
%
%
%
\font\mycmmiten=cmmi10%
\font\mycmmiseven=cmmi7%
\font\mycmmifive=cmmi5%
\newfam\mycmmifam%
\textfont\mycmmifam=\mycmmiten%
\scriptfont\mycmmifam=\mycmmiseven%
\scriptscriptfont\mycmmifam=\mycmmifive%
\def\hexa#1{\ifcase #1 0\or 1\or 2\or 3\or 4\or 5\or 6\or 7\or 8\or 9\or A\or B\or C\or D\or E\or F\fi}%
\mathchardef\mathi="7\hexa\mycmmifam7B%
\mathchardef\mathj="7\hexa\mycmmifam7C%
%
%
\font\mymsbmten=msbm10 at 8pt%
\font\mymsbmseven=msbm7 at 5.6pt
\font\mymsbmfive=msbm5 at 4pt%
\newfam\mymsbmfam%
\textfont\mymsbmfam=\mymsbmten%
\scriptfont\mymsbmfam=\mymsbmseven%
\scriptscriptfont\mymsbmfam=\mymsbmfive%
\mathchardef\mybeth="7\hexa\mymsbmfam69%
\mathchardef\mygimmel="7\hexa\mymsbmfam6A%
\mathchardef\mydaleth="7\hexa\mymsbmfam6B%
%
%
%
%
\def\proof{{\noindent\bf Proof:\ }}%
\def\remark{{\noindent\bf Remark:\ }}%
\def\qed{~$\square$}%
\def\makeop#1{\global\expandafter\def\csname op#1\endcsname{{\text{{\rm #1}}}}}%
\def\makeopsmall#1{\global\expandafter\def\csname op#1\endcsname{{\text{{\rm \lowercase{#1}}}}}}%
%
%
%
%
%
%
\makeop{Ext}%
\makeop{Int}%
\makeop{Dist}%
\makeop{Diam}%
\makeop{Length}%
%
%
%
%
%
%
%
%
\def\msup{\mathop{{\text{Sup}}}}%
%
%
%
\makeop{Dim}%
\makeop{Ker}%
\makeop{Coker}%
\makeop{Tr}%
\makeop{Adj}%
\makeop{Det}%
\makeop{End}%
\makeop{Lin}%
\makeop{Symm}%
\makeop{Mult}%
%
%
\makeop{dx}%
\makeop{dy}%
\makeop{dz}%
\makeop{dt}%
\makeop{dVol}%
\makeop{dArea}%
\makeop{Supp}%
\makeop{Hess}%
\makeop{Lip}%
%
%
\makeop{Re}%
\makeop{Im}%
\makeop{Arg}%
\makeop{Log}%
\makeop{Exp}%
%
%
\makeopsmall{Cos}%
\makeopsmall{Sin}%
\makeopsmall{Tan}%
\makeopsmall{Sec}%
\makeopsmall{Cosec}%
\makeopsmall{Cot}%
\makeopsmall{ArcCos}%
\makeopsmall{ArcSin}%
\makeopsmall{ArcTan}%
\makeopsmall{ArcSec}%
\makeopsmall{ArcCosec}%
\makeopsmall{ArcCot}%
%
%
\makeopsmall{Cosh}%
\makeopsmall{Sinh}%
\makeopsmall{Tanh}%
\makeopsmall{ArcCosh}%
\makeopsmall{ArcSinh}%
\makeopsmall{ArcTanh}%
%
%
\makeop{Vol}%
\makeop{Area}%
\makeop{Riem}%
\makeop{Ric}%
\makeop{Scal}%
\makeop{Euc}%
\makeop{Imm}%
\makeop{Emb}%
%
%
\makeop{Id}%
\makeop{Ad}%
\makeop{O}%
\makeop{SO}%
\makeop{SL}%
\makeop{GL}%
\makeop{Conf}%
\makeop{Homeo}%
\makeop{Diff}%
\makeop{Isom}%
%
%
\makeop{Ind}%
\makeop{Sig}%
\makeop{Spec}%
%
%
\makeop{Conv}%
\makeop{Max}%
\makeop{Min}%
\makeop{Mod}%
\makeop{Deg}%
\makeop{loc}%
%
%
%
%
%
%
%
%
%
%
%
%
%
 %
%
%
%
%
%
\newref{BrendleEichmair}{Brendle S., Eichmair M., Large outlying stable constant mean curvature spheres in initial data sets, {\sl Invent. Math.}, {\bf 197}, (2014), no. 3, 663--682}
\newref{Chavel}{Chavel I., {\sl Riemannian geometry}, Cambridge Studies in Advanced Mathematics, {\bf 98}, Cambridge University Press, Cambridge, (2006)}
\newref{EichmairMetzgerI}{Eichmair M., Metzger J., On large volume preserving stable (CMC) surfaces in initial data sets, {\sl J. Diff. Geom}, {\bf 91}, (2012), no. 1, 81--102}
\newref{EichmairMetzgerII}{Eichmair M., Metzger J., Large isoperimetric surfaces in initial data sets, {\sl J. Diff. Geom.}, {\bf 94}, (2013), no. 1, 159--186}
\newref{EichmairMetzgerIII}{Eichmair M., Metzger J., Unique isoperimetric foliations of asymptotically flat manifolds in all dimensions, {\sl Invent. Math.}, {\bf 194}, (2013), no. 3, 591--630}
\newref{Krylov}{Krylov N. V., {\sl Lectures on elliptic and parabolic equations in Sobolev spaces}, Graduate Studies in Mathematics, {\bf 96}, American Mathematical Society, Providence, RI, (2008)}
\newref{MaxNunSmi}{Maximo D., Nunes I., Smith G., Free boundary minimal annuli in convex three-manifolds, arXiv:1312.5392}
\newref{NardulliI}{Nardulli S., The isoperimetric profile of a smooth riemannian manifold for small volumes, {\sl Ann. Global Anal. Geom.}, {\bf 36}, (2009), no. 2, 111--131}
\newref{NardulliII}{Nardulli S., The isoperimetric profile of a non-compact riemannian manifold for small volumes, {\sl Calc. Var. Partial Differential Equations}, {\bf 49}, (2014), no. 1-2, 173--195}
\newref{Salamon}{Robbin J. W., Salamon D. A., The spectral flow and the Maslov index, {\sl Bull. London Math. Soc.}, {\bf 27}, (1995), 1--33}
\newref{RosSmi}{Rosenberg H., Smith G., Degree theory of immersed hypersurfaces, arXiv:1010.1879}
\newref{Rudin}{Rudin W., {\sl Principles of mathematical analysis}, International Series in Pure and Applied Mathematics, McGraw-Hill Book Co., New York, Auckland, D\"usseldorf, (1976)}
\newref{Schwarz}{Schwarz M., {\sl Morse homology}, Progress in Mathematics, {\bf 111}, Birkh\"auser Verlag, Basel, 1993}
\newref{SmiCCH}{Smith G., Constant curvature hypersurfaces and the Euler characteristic, arXiv:1103.3235}
\newref{SmiEMCFI}{Smith G., Eternal forced mean curvature flows I: a compactness result, {\sl Geom. Dedicata}, {\bf 176}, (2015), 11--29}
\newref{White}{White B., The space of minimal submanifolds for varying riemannian metrics, {\sl Indiana
Math. Journal}, {\bf 40}, (1991), no.1, 161--200}
\newref{Ye}{Ye R., Foliation by constant mean curvature spheres, {\sl Pacific J. Math.}, {\bf 147}, (1991), no. 2, 381--396}
\catcode`\@=11
\def\triplealign#1{\null\,\vcenter{\openup1\jot \m@th %
\ialign{\strut\hfil$\displaystyle{##}$&$\displaystyle{{}##}\hfil$&$\displaystyle{{}##}$\hfil\crcr#1\crcr}}\,}
\def\multiline#1{\null\,\vcenter{\openup1\jot \m@th %
\ialign{\strut$\displaystyle{##}$\hfil&$\displaystyle{{}##}$\hfil\crcr#1\crcr}}\,}
\catcode`\@=12
\def\Pagetitle{\hfil}
\def\Pagefooter{\hfil{\myfontdefault\folio}\hfil}
\null \vfill
\def\centre{\rightskip=0pt plus 1fil \leftskip=0pt plus 1fil \spaceskip=.3333em \xspaceskip=.5em \parfillskip=0em \parindent=0em}%
\def\textmonth#1{\ifcase#1\or January\or Febuary\or March\or April\or May\or June\or July\or August\or September\or October\or November\or December\fi}
\font\abstracttitlefont=cmr10 at 14pt {\abstracttitlefont\centre
Eternal Forced Mean Curvature Flows II - Existence.\par}
\bigskip
{\centre 23rd August 2015\par}
\bigskip
{\centre Graham Smith\par}
\bigskip
{\centre Instituto de Matem\'atica,\par
UFRJ, Av. Athos da Silveira Ramos 149,\par
Centro de Tecnologia - Bloco C,\par
Cidade Universit\'aria - Ilha do Fund\~ao,\par
Caixa Postal 68530, 21941-909,\par
Rio de Janeiro, RJ - BRASIL\par}
\bigskip
\noindent{\bf Abstract:\ }We show that under suitable non-degeneracy conditions, complete gradient flow lines of the scalar curvature functional of a riemannian manifold perturb into eternal forced mean curvature flows with large forcing term.
%
%
\par
\vfill
\nextoddpage
\global\pageno=1
\myfontdefault
\def\Pagetitle{\hfil Eternal Forced Mean Curvature Flows II - Existence\hfil}
\def\Pagefooter{\hfil{\myfontdefault\folio}\hfil}
\makeop{in}
\newhead{Introduction}[Introduction]
\newsubhead{Background}[Background]
In \cite{Ye}, Ye shows how non-degenerate critical points of the scalar curvature function of a riemannian manifold perturb into families of convex embedded spheres of arbitrary large constant mean curvature inside that manifold. While this result has been shown to have significant applications in the study of the isoperimetric problem (c.f., for example, \cite{BrendleEichmair}, \cite{EichmairMetzgerI}, \cite{EichmairMetzgerII}, \cite{EichmairMetzgerIII}, \cite{NardulliI} and \cite{NardulliII}), its applications to the study of the differential topologies of spaces of immmersed and embedded submanifolds have been less exploited. However, in \cite{SmiCCH}, we show how Ye's result implies that - in heuristic terms - the Euler characteristic of the space of convex Alexandrov embedded spheres inside a given manifold is equal to $(-1)$ times the Euler characteristic of that manifold. This has applications to the study of existence, and to some measure, uniqueness, of Alexandrov embedded spheres of constant curvature for many different notions of curvature.
\par
However, if our aim is to prove existence, then the results of \cite{SmiCCH} are unsatisfactory when the Euler characteristic of the ambient manifold vanishes. This happens, for example, when the ambient manifold is $3$-dimensional, which is nonetheless one of the most interesting cases. Furthermore, even when these techniques can be successfully applied to prove existence (as in, for example, \cite{MaxNunSmi}, \cite{RosSmi} or \cite{White}), they still often fall short of optimal results, for there are good topological reasons to believe that - at least generically - there are far more solutions than those whose existence we have managed to prove.
\par
With this in mind, in \cite{SmiEMCFI}, we initiated a programme for the study of the Morse homology of the spaces of immersed and embedded hypersurfaces, where the natural Morse function to be studied is the area functional, or, more generally, the ``Area minus Volume'' functional (defined below), which depends on a parameter $h$, and which we denote by $\Cal{F}_h$. The critical points of $\Cal{F}_h$, which define the chain groups of the Morse complex (c.f. \cite{Schwarz}), are then immersed hypersurfaces of constant mean curvature equal to $h$, and its complete gradient flows, which define the $\partial$ operator of this complex (c.f. \cite{Schwarz}, again), are then eternal forced mean curvature flows with forcing term $h$.
\par
Within this context, Ye's result says that for large values of $h$, non-degenerate critical points of the scalar curvature function map to (in fact, non-degenerate) critical points of $\Cal{F}_h$. In this paper, we prove the corresponding result for complete gradient flows of the scalar curvature function. That is, under suitable non-degeneracy conditions, we show that for sufficiently large values of $h$, these flows map to complete gradient flows of $\Cal{F}_h$. Combined with a suitable converse (that is, a concentration result), which has been proven in Ye's case, but which we have not yet proven here, this would mean that for large values of $h$, the entire Morse complex of the scalar curvature functional maps to the Morse complex of $\Cal{F}_h$. This would make the two isomorphic, thereby yielding an explicit description of the Morse homology of the space of Alexandrov embedded spheres. Since the number of constant mean curvature immersed spheres should be bounded below by the sum of the Betti numbers of this homology, we should thereby obtain stronger existence results for such hypersurfaces than those that are currently known.
\newsubhead{Notation, Terminology and Main Result}[NotationTerminologyAndMainResult]
Let $M:=M^{m+1}$ be a complete $(m+1)$-dimensional riemannian manifold. Let $S$ be its scalar curvature function, where, throughout the paper, we adopt the convention which normalises all curvature functions so that the unit sphere in Euclidean space always has positive unit curvature. Let $\gamma:\Bbb{R}\rightarrow M$ solve the non-linear ODE
$$
\dot{\gamma} + \frac{(m+1)}{2(m+3)}\nabla S = 0,\eqnum{\nexteqnno[GradientFlowEquation]}
$$
so that $\gamma$ is (up to reparametrisation) a {\bf complete gradient flow line} of $S$. Consider the linearisation $L$ of \eqnref{GradientFlowEquation} about $\gamma$. This is a linear ordinary differential operator which maps $\Gamma(\gamma^*TM)$ to itself and is given by
$$
L = \frac{\partial}{\partial t} + \frac{(m+1)}{2(m+3)}\opHess(S).\eqnum{\nexteqnno[LinearisedGradientFlowEquation]}
$$
We now recall that $S$ is said to be of {\bf Morse} type whenever all of its critical points are non-degenerate. In this case, if $\gamma$ has relatively compact image, then $\gamma(t)$ converges towards critical points of $S$ as $t$ tends to $\pm\infty$. Furthermore (c.f. \cite{Salamon}), $L$ defines a Fredholm mapping from the space of H\"older differentiable sections, $\Gamma^{k+1,\alpha}(\gamma^*TM)$, into $\Gamma^{k,\alpha}(\gamma^*TM)$, whose Fredholm index is equal to the difference of the Morse indices of the two end-points of $\gamma$. We then say that $\gamma$ is {\bf non-degenerate} whenever $L$ is surjective, and we say that $S$ is of {\bf Morse-Smale} type whenever, in addition to all of its critical points being non-degenerate, all of its complete gradient flows which have relatively compact image are also non-degenerate. This is the property that we require for the Morse complex of $S$ to be well-defined (c.f. \cite{Schwarz}). There is no shortage of metrics whose scalar curvature function has this property. Indeed, the set of all such metrics is generic (that is, in the second category in the sense of Baire) within any conformal class.
\par
Now let $B^{m+1}$ and $S^m$ be respectively the closed unit ball and the unit sphere in $\Bbb{R}^{m+1}$. Let $\hat{\Cal{E}}$ be the space of smooth immersions of $B^{m+1}$ into $M$ and let $\Cal{E}$ be the quotient of this space under the action of the group of smooth orientation preserving diffeomorphisms of $B^{m+1}$. It is usual to  identify an immersion in $\hat{\Cal{E}}$ with its equivalence class in $\Cal{E}$. By a slight abuse of terminology, for each $e\in\Cal{E}$, we define $\opVol(e)$ and $\opArea(e)$ to be respectively the volumes of $B^{m+1}$ and $S^m$ with respect to the metric $e^*g$. For all $h>0$, we now define the ``Area minus Volume'' functional by
$$
\Cal{F}_h(e) := \opArea(e) - h\opVol(e).\eqnum{\nexteqnno[AreaMinusVolumeFunctional]}
$$
\par
Many properties of the immersion $e$ are actually determined by its restriction to $S^m$. Indeed, the restriction operator actually defines a local homeomorphism from $\Cal{E}$ into the space of reparametrisation equivalence classes of immersions of $S^m$ into $M$, whose image is the space of {\bf Alexandrov embeddings} of $S^m$ into $M$. Furthermore, the embedding $e:B^{m+1}\rightarrow M$ is a critical point of $\Cal{F}_h$ whenever its restriction to $S^m$ has constant mean curvature equal to $h$. Likewise, the family $e:\Bbb{R}\times B^{m+1}\rightarrow M$ is an $L^2$ gradient flow of $\Cal{F}_h$ whenever its restriction to $\Bbb{R}\times S^m$ is a forced mean curvature flow with forcing term $h$. That is, whenever this restriction satisfies
$$
\left\langle\frac{\partial}{\partial t}e,N_t\right\rangle + H_t - h = 0,\eqnum{\nexteqnno[ForcedMeanCurvatureFlowEquation]}
$$
where $N_t$ and $H_t$ are respectively the outward-pointing unit normal vector field and the mean curvature of the restriction of $e_t:=e(t,\cdot)$ to $\Bbb{R}\times S^m$.
\par
We now introduce the mechanism by which complete gradient flow lines of $S$ perturb to eternal forced mean curvature flows. Let $\gamma$ be a complete gradient flow line of $S$. Using parallel transport, we identify the bundle $\gamma^*TM$ with the trivial bundle $\Bbb{R}\times\Bbb{R}^{m+1}$, and we define $\opExp:\Bbb{R}\times\Bbb{R}^{m+1}\rightarrow M$ such that, for all $t$, $\opExp_t:=\opExp(t,\cdot)$ is the exponential map of $M$ about the point $\gamma(t)$. Now, following \cite{Ye}, for all $s>0$, for all $Y:\Bbb{R}\rightarrow\Bbb{R}^{m+1}$ and for all $f:\Bbb{R}\times S^m\rightarrow]0,\infty[$, we define the function $e(s,Y,f):\Bbb{R}\times S^m\rightarrow M$ by
$$
e(s,Y,f)(t,x) = \opExp_t(sY(t) + s(1 + s^2f(t,x))x).\eqnum{\nexteqnno[PerturbedFlowOfSpheres]}
$$
Heuristically, $e(s,Y,f)$ is a smooth family of immersed spheres in $M$ whose centres move along $\gamma$ with a small displacement given by $Y$.
\proclaim{Theorem \nextprocno}
\noindent If $S$ is of Morse-Smale type, and if $\gamma$ is a complete gradient flow line of $S$ with relatively compact image, then, for all sufficiently small $s$, there exist $Y:\Bbb{R}\rightarrow\Bbb{R}^{m+1}$ and $f:\Bbb{R}\times S^m\rightarrow]0,\infty[$ such that, up to reparametrisation in time, $e(s,Y,f)$ is an eternal forced mean curvature flow with forcing term $1/s$.
\endproclaim
\proclabel{ThmMainTheorem}
\remark A detailed formal statement of Theorem \procref{ThmMainTheorem} is given in Theorem \procref{ThmExistence} below. In particular, not only do we obtain H\"older estimates for the pair $(Y,f)$, but we also describe in Theorem \procref{ExistenceOfApproximateSolutions}, below, an iterative process for determining asymptotic expansions of these solutions up to arbitrary order.
\newsubhead{Discussion}[Discussion]
Like Ye's result, Theorem \procref{ThmMainTheorem} is proven by first determining formal solutions in the form of asymptotic series, and then perturbing suitably high order partial sums of these series to yield exact solutions. There are, nonetheless, considerable differences between Theorem \procref{ThmMainTheorem} and Ye's result, primarily because Theorem \procref{ThmMainTheorem} is a parabolic, and not an elliptic, problem. Now, on the one hand, since parabolic and elliptic operators are all hypoelliptic, the analytic tools that we use are barely different. However, on the other hand, the time-dependence introduces new - and rather confusing - phenomena as the scale parameter, $s$, tends to zero.
\par
This is perhaps best illustrated by considering the first approximation $Y=0$ and $f=0$. Here, the mean curvature of the sphere $e(s,0,0)(t,\cdot)$ is equal to $1/s + O(s)$, so that the forced mean curvature flow with forcing term $1/s$ should move along the curve $\gamma$ with speed approximately $s$, which trivially tends to $0$. It is perhaps surprising that this scale dependence does not actually introduce any singularities as $s$ tends to $0$. However, a deeper study of the equations involved reveals the role played by operator
$$
Q_s := s^4\frac{\partial}{\partial t} + \frac{1}{m}(m+\overline{\Delta}),\eqnum{\nexteqnno[ScaleDependenceOperator]}
$$
where $\overline{\Delta}$ is the standard Laplacian of the sphere $S^m$. Here the time dependence introduces a fourth power of $s$, and this does affect us in three different ways.
\par
First, Theorem \procref{ThmMainTheorem} becomes a genuine singular perturbation problem. In actual fact, Ye's result, although presented as a singular perturbation problem, transforms, after removing the first few terms and then dividing by a suitable factor, into a regular perturbation problem, which is then directly solved by the inverse function theorem. In the present case, however, when $s=0$, the operator $Q_s$ is no longer hypo-elliptic, and the same simplification no longer applies.
\par
Second, since the Green's operator of $Q_s$ depends on $s$, the terms in the asymptotic series of the formal solution (determined in Theorem \procref{ExistenceOfApproximateSolutions}, below) actually also depend on $s$, so that more care is required in ensuring H\"older bounds which are independent of $s$.
\par
Third, the appropriate functional analytic framework for studying parabolic operators is that of inhomogeneous spaces (introduced here in Section \subheadref{ParabolicOperatorsII}, below). Furthermore, the $s$ dependence of $Q_s$ requires the use of weighted spaces (also defined in Section \subheadref{ParabolicOperatorsII}, below), where what appears to be the most appropriate weighting is in fact slightly counter-intuitive (c.f. the remarks following Lemma \procref{InfiniteDimensionalFactor}).
\par
Finally, in order to develop a Morse homology theory for the space of convex, Alexandrov embedded spheres, two further results are still required. Indeed, it would be necessary to show, first  that the eternal flows obtained here are non-degenerate, and second, that for sufficiently large values of the forcing term, they are the only ones. However, we believe at this stage that it is more interesting to develop a more satisfactory compactness result than that obtained in \cite{SmiEMCFI}, and for this reason we postpone this study to later work.
\newsubhead{Overview of Paper and Acknowledgements}[OverviewAcknowledgements]
This paper is structured as follows. In Section $2$, we develop a formalism for the succinct description of the Taylor series of various well-known geometric functions, and in Section $3$, we extend this formalism in order to describe the functions used in the proof of Theorem \procref{ThmMainTheorem}. Our objective here is to understand the general terms of these series without having to resort to explicit calculations, and, for the sake of completeness, we have studied this problem far more deeply than is actually necessary for our current applications. In Section $4$, we then reformulate these results in the language of asymptotic series. In particular, since the operation of composition by smooth functions yields smooth functionals between H\"older spaces, this immediately yields norm estimates for the functionals of interest to us without any further effort being required.
\par
Having determined the asymptotic expansion of the forced mean curvature flow operator, the rest of Section 4 is devoted to constructing formal solutions and then perturbing these formal solutions into exact solutions. It is here that we introduce the required functional analytic framework, based on the Fredholm theory of parabolic operators over weighted inhomogeneous H\"older spaces (c.f. \cite{Krylov}). In addition, using the theory of spherical harmonics, we improve our norm estimates for every term in the aysmptotic series of our formal solutions. Although this is not strictly necessary, we believe it makes our reasoning a lot cleaner. Finally, once formal solutions have been constructed, a straightforward application of the inverse function theorem yields the desired result.
\par
The results of this paper were established when the author was benefitting from a Marie Curie postdoctoral fellowship in the Centre de Recerca Matem\`atica, Barcelona. The author is grateful to Andrew Clarke for helpful conversations.
\newhead{The Taylor Series of Geometric Functions}[TheTaylorSeriesOfGeometricFunctions]
\newsubhead{Curvature Tensors}[CurvatureTensors]
Let $\Omega$ be the unit ball in $\Bbb{R}^{m+1}$. Let $g$ be a smooth metric over $\Omega$ with Levi-Civita covariant derivative $\nabla$ and Riemann curvature tensor $R$. We suppose that
$$
\nabla_{\partial_r}\partial_r = 0,\ \text{and}\ g(\partial_r,\partial_r) = 1,\eqnum{\nexteqnno[MetricIsExponential]}
$$
where $\partial_r$ here denotes the unit radial vector field. This simply means that $(\Omega,g)$ is an exponential chart of some riemannian manifold. Now denote $\delta_{ij}:=g(0)_{ij}$ and let $\delta^{ij}$ be its metric dual, so that, by \eqnref{MetricIsExponential}, $\delta_{ij}$ is simply the standard euclidean metric over $\Bbb{R}^{m+1}$. Finally, for convenience, we suppose that $\Omega$ is convex in the sense that for all $x,y\in\Omega$, there exists a unique geodesic in $\Omega$ from $x$ to $y$.
\par
We say that a function defined over $\Omega$ is {\bf geometric} when it only depends on the metric $g$. We are interested in the Taylor series about $0$ of such functions, and, in particular, how their coefficients depend on the Riemann curvature tensor. In order to describe this dependence, we introduce the following algebraic formalism. Consider the set of formal tensors $X:=\left\{(R_{i_1i_2i_3}{}^j{}_{;i_4...i_{k+3}})_{k\in\Bbb{N}}\right\}$ where the subscript $;$ here denotes formal covariant differentiation. Observe that all elements of $X$ are covariant of order $1$ and contravariant of order at least $3$. Given two formal tensors, $\rho^a_{bi_1...i_p}$ and $\rho^b_{j_1...j_q}$, which are both covariant of order $1$, define their {\bf matrix product} by $\rho^a_{bi_1...i_p}\rho^b_{j_1...j_q}$, and observe that this product is also covariant of order $1$. Now let $\Cal{R}$ be the vector space with basis the set of all finite formal combinations of elements of $X$ obtained by permutation of indices and matrix multiplication. We call $\Cal{R}$ the space of {\bf curvature tensors}.
\par
For all $k\in\Bbb{Z}$, let $\Cal{R}^k$ be the subspace of $\Cal{R}$ consisting of those elements which are contravariant of order $k+1$. When $\rho\in\Cal{R}^k$, we say that it has {\bf order-difference} $k$. Observe that order-difference is preserved by permutation of indices, and if $\rho$ and $\rho'$ have order-differences $k$ and $k'$ respectively, then their matrix product has order-difference $k+k'$. In particular, since every generator of $\Cal{R}$ has order-difference at least $2$, it follows that for $k<2$, $\Cal{R}^k$ is trivial, and for $k\geq 2$, it is spanned by matrix products of those generators which have order-difference at most $k$. Considerations such as these make it relatively straightforward to determine $\Cal{R}^k$ for all $k$. For example,
$$\eqalign{
\Cal{R}^2 \hfill&
=\langle(R_{i_{\sigma(1)}i_{\sigma(2)}i_{\sigma(3)}}{}^j)_{\sigma\in\Sigma_3}\rangle,\hfill\cr
\Cal{R}^3 \hfill&
=\langle (R_{i_{\sigma(1)}i_{\sigma(2)}i_{\sigma(3)}}{}^j{}_{;i_{\sigma(4)}})_{\sigma\in\Sigma_4}\rangle,\hfill\cr
\Cal{R}^4 \hfill&
=\langle (R_{i_{\sigma(1)}i_{\sigma(2)}i_{\sigma(3)}}{}^j{}_{;i_{\sigma(4)}i_{\sigma(5)}},
R_{pi_{\sigma(1)}i_{\sigma(2)}}{}^jR_{i_{\sigma(3)}i_{\sigma(4)}i_{\sigma(5)}}{}^p,\hfill\cr
&\qquad\qquad R_{i_{\sigma(1)}i_{\sigma(2)}p}{}^jR_{i_{\sigma(3)}i_{\sigma(4)}i_{\sigma(5)}}{}^p)_{\sigma\in\Sigma_5}\rangle,\hfill\cr}
\eqnum{\nexteqnno[SpacesOfCurvatureTensors]}
$$
and so on, where, for all $k$, $\Sigma_k$ denotes the group of permutations of the set $\left\{1,...,k\right\}$.
\par
Identifying elements of $\Cal{R}$ via the symmetries of the Riemann curvature tensor, we obtain
\proclaim{Proposition \nextprocno}
\noindent $\Cal{R}$ is self-adjoint with respect to $\delta$ in the sense that if $\rho^i_{j_1...j_k}$ is an element of $\Cal{R}$, then $\delta^{ia}\delta_{j_lb}\rho^b_{j_1...j_{l-1}aj_{l+1}...j_k}$ identifies with a unique element of $\Cal{R}$ for all $1\leq l\leq k$.
\endproclaim
\proclabel{TheSpaceOfCurvaturePolynomialsIsSelfAdjoint}
\proof It suffices to prove the result for each generator of $\Cal{R}$. We thus show that $\delta^{ia}\delta_{j_lb}R_{j_1j_2j_3}{}^b{}_{;j_4...j_{l-1}aj_{l+1}...j_{k+3}}$ identifies with a unique element of $\Cal{R}$ for all $k$ and for all $1\leq l\leq k+3$. We achieve this by induction on $k$. Indeed, for $k=0$, the result follows directly from the symmetries of the Riemann curvature tensor. For $k=1$, it follows from these symmetries together with the second Bianchi identity. Now suppose that $k\geq 2$. Since the set of generators of $\Cal{R}$ is closed under formal covariant differentiation, so too is $\Cal{R}$, and we may therefore suppose that $l=k+3$. However,
$$\eqalign{
R_{j_1j_2j_3}{}^i{}_{;j_4...j_{k+2}j_l} \hfill&= R_{j_1j_2j_3}{}^i{}_{;j_4...j_{k+1}j_lj_{k+2}} + R_{j_{k+2}j_l p}{}^iR_{j_1j_2j_3}{}^p{}_{;j_4...j_{k+1}}\hfill\cr
&\qquad-\sum_{b=1}^{k+1}R_{j_{k+2}j_lj_b}{}^aR_{j_1j_2j_3}{}^i{}_{;j_4...j_{b-1}aj_{b+1}...j_{k+1}},\hfill\cr}
$$
and the result now follows by induction.\qed
\medskip
\noindent The significance of Proposition \procref{TheSpaceOfCurvaturePolynomialsIsSelfAdjoint} lies in the fact that although geometric functions are defined in terms of the metric, they can be approximated purely in terms of curvature tensors, as we will see presently.
\par
Finally, denote
$$
\overline{\Cal{R}} := \Cal{R}\oplus\langle\delta^i_j\rangle,\eqnum{\nexteqnno[DefnOfOverlineR]}
$$
and, for all $k$, define $\overline{\Cal{R}}^k$ as before. We also call elements of $\overline{\Cal{R}}$ {\bf curvature tensors}. Observe that $\overline{\Cal{R}}$ is also closed under matrix multiplication. Furthermore, $\overline{\Cal{R}}^0=\langle\delta^i_j\rangle$, and, for all $k\neq 0$, $\overline{\Cal{R}}^k=\Cal{R}^k$.
\newsubhead{Curvature Polynomials}[CurvaturePolynomials]
Let $\underline{X}:=(X_1,...,X_n)$ be a vector of formal variables each taking values in $\Bbb{R}^{m+1}$. For $\rho\in\Cal{R}^k$ and for $0\leq r_1+...+r_n\leq k+1$, define the formal polynomial
$$
(\rho_{r_1,...,r_n})^i_{j_{(r_1+...+r_n)+1}...j_{k+1}} := \rho^i_{j_1...j_{k+1}}X_1^{j_1}...X_1^{j_{r_1}}...X_n^{j_{(r_1+...+r_{n-1}) + 1}}...X_n^{(r_1+...+r_n)},
\eqnum{\nexteqnno[FundamentalCurvaturePolynomial]}
$$
where $X^i_j$ denotes the $i$'th component of the vector $X_j$. Abusing notation, let $\Cal{R}[\underline{X}]$ be the vector space with basis the set of all such formal polynomials. We call $\Cal{R}[\underline{X}]$ the space of {\bf curvature polynomials}. Observe that $\Cal{R}[\underline{X}]$ is closed under matrix multiplication, although it is not always possible to multiply two given elements (indeed, two elements which are both covariant of order $1$ and contravariant of order $0$ cannot be multiplied). Furthermore, since $\Cal{R}$ is self-adjoint with respect to $\delta$, so too is $\Cal{R}[\underline{X}]$ in the sense that if $P^{i}_{j_1...j_k}$ is an element of $\Cal{R}[\underline{X}]$, then $\delta^{ia}\delta_{j_lb}P^b_{j_1...j_{l-1}aj_{l+1}...j_k}$ identifies with a unique element of $\Cal{R}[\underline{X}]$ for all $1\leq l\leq k$.
\par
For $k\in\Bbb{Z}$ and for $\underline{r}:=(r_1,...,r_n)\in\Bbb{N}^n$, let $\Cal{R}^k_{\underline{r}}[\underline{X}]$ denote the subspace of $\Cal{R}[\underline{X}]$ consisting of those elements which are contravariant of order $k+1$ and homogeneous of degree $r_i$ in $X_i$ for each $i$. Likewise, denote
$$
\Cal{R}^k[\underline{X}] := \oplus_{\underline{r}}\Cal{R}^k_{\underline{r}}[\underline{X}].
\eqnum{\nexteqnno[SpaceOfCurvaturePolynomials]}
$$
When $P\in\Cal{R}^k_{\underline{r}}[\underline{X}]$, we say that it has {\bf order-difference} $k$ and {\bf degree} $\underline{r}$. As before, permutation of indices preserves order-difference, and if $P$ and $P'$ have order-differences $k$ and $k'$ respectively then their matrix product has order-difference $k+k'$.
\par
Throughout most of this section, we will only be concerned with the case where $n=1$. Here we have
\proclaim{Proposition \nextprocno}
\noindent If $r>k$ and if $\rho\in\Cal{R}^k$, then $\rho_r=0$. In particular, $\Cal{R}^k[X]$ is non-trivial only if $k\geq 0$.
\endproclaim
\proclabel{CurvaturePolynomialsInOneVariableHaveNonNegativeOrderDifference}
\proof It suffices to prove the result when $\rho$ is a generator of $\Cal{R}$. However, for each $k$, by symmetry, $R_{j_1j_2j_3}{}^i{}_{;j_4...j_k}X^{j_1}...X^{j_k}=0$, and the result follows.\qed
\medskip
\noindent Proposition \procref{CurvaturePolynomialsInOneVariableHaveNonNegativeOrderDifference} implies that every element of $\Cal{R}^0[X]$ is a finite sum of matrix products of those generators of $\Cal{R}[X]$ which are of order-difference $0$, that is, formal polynomials of the form $R_{p_1ip_2}{}^j{}_{;p_3...p_{k+2}}X^{p_1}...X^{p_{k+2}}$, where $k$ varies over all non-negative integers. By considerations such as these, we obtain, for example,
$$\eqalign{
\Cal{R}^0_0[X] \hfill&= 0,\hfill\cr
\Cal{R}^0_1[X] \hfill&= 0,\hfill\cr
\Cal{R}^0_2[X] \hfill&= \langle R_{piq}{}^jX^pX^q \rangle,\hfill\cr
\Cal{R}^0_3[X] \hfill&= \langle R_{piq}{}^j{}_{;r}X^pX^qX^r \rangle,\hfill\cr}
\eqnum{\nexteqnno[ExamplesOfSpacesOfCurvaturePolynomials]}
$$
and so on. Likewise, every element of $\Cal{R}^1[X]$ is a finite sum of matrix products of generators all but one of which are elements of $\Cal{R}^0[X]$ and the remaining one of which is an element of $\Cal{R}^1[X]$, and we obtain,
$$\eqalign{
\Cal{R}^1_0[X] \hfill&= 0,\hfill\cr
\Cal{R}^1_1[X] \hfill&= \langle (R_{pi_{\sigma(1)}i_{\sigma(2)}}{}^j X^p, R_{i_{\sigma(1)}i_{\sigma(2)}p}{}^jX^p)_{\sigma\in\Sigma_2}\rangle,\hfill\cr
\Cal{R}^1_2[X] \hfill&= \langle (R_{p i_{\sigma(1)} q}{}^j{}_{;i_{\sigma(2)}}X^pX^q, R_{pi_{\sigma(1)}i_{\sigma(2)}}{}^j{}_{;q}X^pX^q, R_{i_{\sigma(1)}i_{\sigma(2)}p}{}^j{}_{;q}X^pX^q)_{\sigma\in\Sigma_2}\rangle,\hfill\cr}
\eqnum{\nexteqnno[MoreExamplesOfSpacesOfCurvaturePolynomials]}
$$
and so on. In summary, it is relatively straightforward to determine $\Cal{R}^k_r[X]$ for all $k$ and for all $r$.
\par
For general $n$, since $\Cal{R}^k$ is trivial for $k<2$, $\Cal{R}^{-1}_{\underline{r}}[\underline{X}]$ is trivial for $r_1+...+r_n\leq 2$ and $\Cal{R}^0_{\underline{r}}[\underline{X}]$ is trivial for $r_1+...+r_n\leq 1$. This observation will play an important role in the sequel.
\par
Finally, as before, denote
$$
\overline{\Cal{R}}[\underline{X}] = \Cal{R}[\underline{X}]\oplus\langle\delta^i_j\rangle\oplus\langle X^i_j\rangle,\eqnum{\nexteqnno[ExtendedCurvaturePolynomials]}
$$
where $X^i_j$ denotes the $i$'th component of the vector $X_j$. For all $k$, and for all $\underline{r}$, define $\overline{\Cal{R}}^k_{\underline{r}}[\underline{X}]$ as before. We also call elements of $\overline{\Cal{R}}[\underline{X}]$ {\bf curvature polynomials}. Observe that $\overline{\Cal{R}}[\underline{X}]$ is also closed under matrix multiplication. Furthermore,
$$\eqalign{
\overline{\Cal{R}}^{-1}[\underline{X}] &= \Cal{R}^{-1}[\underline{X}]\oplus\langle X_1,...,X_n\rangle,\hfill\cr
\overline{\Cal{R}}^0[\underline{X}] &= \Cal{R}^0[\underline{X}]\oplus\langle\delta^i_j\rangle,\hfill\cr}
\eqnum{\nexteqnno[ExamplesOfSpacesOfExtendedCurvaturePolynomials]}
$$
and $\overline{\Cal{R}}^k[\underline{X}]=\Cal{R}^k[\underline{X}]$ for all other values of $k$.
\newsubhead{General Properties of Taylor Series}[GeneralPropertiesOfTaylorSeries]
As before, let $\underline{X}:=(X_1,...,X_n)$ be a vector of formal variables taking values in $\Bbb{R}^{m+1}$. Abusing notation, let $A[\underline{X}]$ be an algebra of formal polynomials in $\underline{X}$, and let $A[[\underline{X}]]$ be the algebra of formal power series in $\underline{X}$ all of whose partial sums are elements of $A[\underline{X}]$. For such a formal power series, $F$, and for every non-negative integer, $k$, denote by $[F]_k$ its partial sum of order $k$.
\par
Recall that for all real $\alpha$, the binomial theorem furnishes a sequence $(a_{k,\alpha})$ of real numbers such that for $x\in]-1,1[$,
$$
(1+x)^\alpha = \sum_{k=0}^\infty a_{k,\alpha}x^k.
\eqnum{\nexteqnno[BinomialTheorem]}
$$
Consequently, if the algebra $A[\underline{X}]$ contains an identity, which we always denote by $I$, then for all formal power series $F$ in $A[[\underline{X}]]$ with $F(0)=I$, and for any real exponent, $\alpha$, we define
$$
F^\alpha := \sum_{k=0}^\infty a_{k,\alpha}(F - I)^k.
\eqnum{\nexteqnno[DefinitionOfPowerOfSeriesByBinomialTheorem]}
$$
\proclaim{Proposition \nextprocno}
\noindent Let $F$ be a formal power series in $X$. If $F$ belongs to $A[[\underline{X}]]$, and if $F(0)=I$, then $F^\alpha$ also belongs to $A[[\underline{X}]]$ for all real $\alpha$.
\endproclaim
\proclabel{PowersOfCurvatureSeriesAreAlsoCurvatureSeries}
\proof Denote $G:=F-I$. For all $k$, $G^k\in A[[\underline{X}]]$ and since $G(0)=0$, $[G^k]_l=0$ for all $l<k$. Thus, for all $\alpha$ and for all $l$,
$$
[F^\alpha]_l = \left[\sum_{k=0}^\infty a_{k,\alpha}G^k\right]_l = \sum_{k=0}^l a_{k,\alpha}[G^k]_l \in A[\underline{X}],
$$
and so $F^\alpha\in A[[\underline{X}]]$, as desired.\qed
\medskip
Now let $\underline{T}:=(T_1,...,T_n)$ be a vector of formal variables taking values in $\Bbb{R}$, and let $A[\underline{X}][[\underline{T}]]$ denote the algebra of formal power series in $\underline{T}$ all of whose coefficients are elements of $A[\underline{X}]$.
\proclaim{Proposition \nextprocno}
\noindent Let $F$ be a formal power series in $\underline{X}$, and define $G(\underline{X},\underline{T}):=F(T_1X_1,....,T_nX_n)$. If $G$ belongs to $A[\underline{X}][[\underline{T}]]$, then $F$ belongs to $A[[\underline{X}]]$.
\endproclaim
\proclabel{TestForBeingCurvatureSeries}
\proof By hypothesis,
$$
G = \sum_{\underline{k}}\frac{1}{k_1!...k_n!}T_1^{k_1}...T_n^{k_n}P_{\underline{k}}(\underline{X}),
$$
where, for all $\underline{k}$, the formal polynomial $P_{\underline{k}}$ belongs to $A[\underline{X}]$. Now consider the formal derivatives of $F$ and $G$ with respect to $\underline{X}$ and $\underline{T}$ respectively. By the chain rule,
$$
\frac{\partial^{k_1}...\partial^{k_n}F}{\partial X_1^{k_1}...\partial X_n^{k_n}}(0)(X_1^{\otimes k_1}\otimes...\otimes X_m^{\otimes k_m}) = \frac{\partial^{k_1}...\partial^{k_n}G}{\partial T_1^{k_1}...\partial T_n^{k_n}}(0,\underline{X}) = P_{\underline{k}}(\underline{X}) \in A[\underline{X}].
$$
It follows that every partial sum of $F$ belongs to $A[\underline{X}]$, and so $F$ belongs to $A[[\underline{X}]]$, as desired.\qed
\newsubhead{Tensor-Valued Geometric Functions}[TensorValuedGeometricFunction]
For all $p,q\in\Bbb{N}$, let $T^{p,q}:=T^{p,q}(\Bbb{R}^{m+1})$ be the space of tensors over $\Bbb{R}^{m+1}$ which are covariant of order $p$ and contravariant of order $q$. Consider a function $f:\Omega\rightarrow T^{1,k+1}$, and denote by $[f]$ its Taylor series. In the present context, the statement that $[f]$ belongs to $R^k[[X]]$ means that the Taylor series of $[f]$ about $0$ is given by
$$
f(x) \sim \sum_{r=0}^\infty R_r(x),
$$
where, for all $r$, $R_r$ is a curvature polynomial of order-difference $k$ and degree $r$.
\par
Now observe that $T^{1,1}$ naturally identifies with $\opEnd(\Bbb{R}^{m+1})$. In particular, since matrix multiplication coincides with the usual notion of matrix multiplication in this case, the space $\overline{\Cal{R}}^0[\underline{X}]$ is also closed with respect to this product, and therefore constitutes an algebra.
\par
Let $M:\Omega\rightarrow\opEnd(\Bbb{R}^{m+1})$ be such that for all $x\in\Omega$ and for every vector $U$, $M(x)U$ is the parallel transport of $U$ along the radial line from $x$ to $0$. Classical Jacobi field techniques (c.f. \cite{Chavel}) readily yield
$$
M^i_j(x) \sim \delta^i_j + \frac{1}{2}R_{pjq}{}^ix^px^q + \frac{1}{3}R_{pjq}{}^i{}_{;r}x^px^qx^r + O(x^4).\eqnum{\nexteqnno[EqnParallelTransport]}
$$
More generally,
\noskipproclaim{Proposition \nextprocno}
$$
[M]\in\overline{\Cal{R}}^0[[X]].\eqnum{\nexteqnno[TaylorSeriesOfM]}
$$
\endproclaim
\proclabel{TaylorSeriesOfM}
\proof Fix a point $x_0\in\Omega$ and a vector $U_0\in\Bbb{R}^{m+1}$. Let $x(t):=tx_0$ and $U(t):=tU_0$, so that $x$ is a geodesic and $U$ is a Jacobi field over $x$. We use a dot to denote differentiation and covariant differentiation in the radial direction. By definition, $U(0)=0$, and since $(\Omega,g)$ is an exponential chart, $\dot{U}(0)=U_0$. Furthermore, by the Jacobi field equation, $\ddot{U} = R_{\dot{x}U}\dot{x}$. We now claim that there exist sequences $(P_k)$ and $(Q_k)$ of polynomials over $\opEnd(\Bbb{R}^{n+1})$ such that, for all $k$,
$$
\nabla^{k+2}_{\dot{x}}U = P_k(R(x)(\dot{x}),...,\nabla^kR(x)(\dot{x}))U + Q_k(R(x)(\dot{x}),...,\nabla^{k-1}R(x)(\dot{x}))\dot{U},
$$
where, for all $l$,
$$
\nabla^lR(x)(\dot{x}) := R(x)_{p_1jp_2}{}^i{}_{;p_3...p_{l+2}}\dot{x}^{p_1}...\dot{x}^{p_{l+2}}.
$$
This holds for $k=0$ by the Jacobi field equation. For $k\geq 0$, using the inductive hypothesis and the fact that $\nabla_{\dot{x}}\dot{x}=0$, we obtain,
$$\eqalign{
\nabla_{\dot{x}}^{k+3}U \hfill
&=\nabla_{\dot{x}}\left(P_k(R(x)(\dot{x}),...,\nabla^kR(x)(\dot{x}))U + Q_k(R(x)(\dot{x}),...,\nabla^{k-1}R(x)(\dot{x}))\dot{U}\right)\hfill\cr
&=P_k(R(x)(\dot{x}),...,\nabla^kR(x)(\dot{x}))\dot{U} + Q_k(R(x)(\dot{x}),...,\nabla^{k-1}R(x)(\dot{x}))\ddot{U}\hfill\cr
&\qquad + \sum_{l=0}^kP_{k,l}(R(x)(\dot{x}),...,\nabla^kR(x)(\dot{x}),\nabla^{l+1}R(x)(\dot{x}))U\hfill\cr
&\qquad + \sum_{l=0}^{k-1}Q_{k,l}(R(x)(\dot{x}),...,\nabla^{k-1}R(x)(\dot{x}),\nabla^{l+1}R(x)(\dot{x}))\dot{U},\hfill\cr
}$$
for suitable sequences of polynomials $(P_{k,l})$ and $(Q_{k,l})$. However, by the Jacobi field equation again, $\ddot{U}=R_{\dot{x}U}\dot{x}$, and the assertion follows by induction. Observe, furthermore, that for all $k$, the zeroeth order terms of $P_k$ and $Q_k$ both vanish. Substituting $t=0$ now yields
$$
(\nabla^{k+2}_{x_0}U)(0) = Q_k(R(0)(x_0),...,\nabla^{k-1}R(0)(x_0))U_0.
$$
However, for all $k$,
$$
\partial_t^k tM(tx_0)U_0|_{t=0} = (\nabla_{x_0}^kU)(0),
$$
so that, by Taylor's Theorem,
$$
[M](TX) = \opId + \sum_{k=2}^\infty\frac{T^k}{(k+1)!}Q_{k-1}(R(0)(X),...,\nabla^{k-2}R(0)(X))\in \overline{\Cal{R}}[X][[T]],
$$
and the result now follows by Proposition \procref{TestForBeingCurvatureSeries}.\qed
\medskip
Let $A,B:\Omega\rightarrow\opEnd(\Bbb{R}^{m+1})$ be such that, for all $x$,
$$\triplealign{
g_{ij}(x) &= A^p_i(x)\delta_{pj} &= \delta_{ip}A^p_j(x),\ \&\cr
g^{ij}(x) &= B^i_p(x)\delta^{pj} &= \delta^{ip}B^j_p(x),\cr}\eqnum{\nexteqnno[EqnMetric]}
$$
where $g^{ij}(x)$ here denotes the metric inverse of $g_{ij}(x)$. Using the same Jacobi field techniques as before, we obtain
$$\eqalign{
A^i_j(x) &\sim \delta^i_j + \frac{1}{3}R_{pjq}{}^ix^px^q + \frac{1}{6}R_{pjq}{}^i{}_{;r}x^px^qx^r + O(x^4),\ \&\cr
B^i_j(x) &\sim \delta^i_j - \frac{1}{3}R_{pjq}{}^ix^px^q - \frac{1}{6}R_{pjq}{}^i{}_{;r}x^px^qx^r + O(x^4).\cr}\eqnum{\nexteqnno[EqnShortAsymptoticSeriesForAAndB]}
$$
More generally,
\noskipproclaim{Proposition \nextprocno}
$$
[A],[B]\in\overline{\Cal{R}}^0[[X]].\eqnum{\nexteqnno[TaylorSeriesOfAAndB]}
$$
\endproclaim
\proclabel{TaylorSeriesOfAAndB}
\proof For every point $x$ in $\Omega$ and for all vectors $U$ and $V$ in $\Bbb{R}^{m+1}$,
$$
g(x)(U,V) = g(0)(M(x)U,M(x)V) = \langle M(x)U,M(x)V\rangle = \langle M^*(x)M(x)U,V\rangle.
$$
Since $U$ and $V$ are arbitrary, it follows that $A=M^*M$. However, since $\overline{\Cal{R}}[X]$ is self-adjoint with respect to $\delta$, $[M^*]$ belongs to $\overline{\Cal{R}}[[X]]$ and therefore so too does $[A]=[M^*][M]$. Finally, since $[A](0)=A(0)=I$, by Proposition \procref{PowersOfCurvatureSeriesAreAlsoCurvatureSeries}, $[B]=[A^{-1}]=[A]^{-1}$ also belongs to $\overline{\Cal{R}}[[X]]$, and this completes the proof.\qed
\medskip
Let $\Gamma:\Omega\rightarrow T^{1,2}$ be the Christoffel symbol of the Levi-Civita covariant derivative of $g$. That is,
$$
\Gamma^k_{ij}(x)\partial_k := \nabla_{\partial_i}\partial_j - D_{\partial_i}\partial_j,\eqnum{\nexteqnno[EqnChristoffelSymbol]}
$$
where $D$ denotes the canonical differentiation operator over $\Bbb{R}^{m+1}$. Recall that $\Gamma$ is symmetric in $i$ and $j$. Furthermore, using the same Jacobi field techniques once again, we obtain
$$
\Gamma^k_{ii}(x) \sim \frac{2}{3}R_{pii}{}^kx_p + \frac{5}{12}R_{pii}{}^k{}_{;q}x^px^q + \frac{1}{12}R_{piq}{}^k{}_{;i}x^px^q + O(x^3).\eqnum{\nexteqnno[EqnShortSeriesForChristoffelSymbol]}
$$
More generally,
\noskipproclaim{Proposition \nextprocno}
$$
[\Gamma]\in\overline{\Cal{R}}^1[[X]].\eqnum{\nexteqnno[TaylorSeriesOfChristoffelSymbol]}
$$
\endproclaim
\proclabel{TaylorSeriesOfChristoffelSymbol}
\proof By the Koszul formula, for all vectors $U$, $V$ and $W$ in $\Bbb{R}^{m+1}$ and for every point $x$ in $\Omega$,
$$
2\langle A(x)\Gamma(x)(U,V),W\rangle = \langle DA(x;U)V,W\rangle + \langle DA(x;V)U,W\rangle - \langle DA(x;W)U,V\rangle.\eqnum{\nexteqnno[KoszulFormula]}
$$
Since $[A]$ belongs to $\overline{\Cal{R}}[[X]]$, its formal derivative, $D[A]=[DA]$ also belongs to $\overline{\Cal{R}}[[X]]$. Now let $\Phi:\Omega\rightarrow T^{1,2}$ be such that
$$
\langle\Phi(x)(U,V),W\rangle = \langle DA(x;V)U,W\rangle.
$$
Since $\overline{\Cal{R}}[[X]]$ is self-adjoint with respect to $\delta$, $[\Phi]$ also belongs to $\overline{\Cal{R}}[[X]]$, and therefore, by linearity, so too does $[A\Gamma]$. It follows that $[\Gamma]=[B][A][\Gamma]=[B][A\Gamma]$ belongs to $\overline{\Cal{R}}[[X]]$, and this completes the proof.\qed
\newsubhead{The Exponential Map and Parallel Transport}[TheExponentialMapAndParallelTransport]
\noindent Define $\Omega_2\subseteq\Bbb{R}^{m+1}\times\Bbb{R}^{m+1}$ by
$$
\Omega_2 := \left\{(x,y)\ |\ \|x\|+\|y\| < 1\right\}.\eqnum{\nexteqnno[DomainOfExponentialMap]}
$$
Let $\opExp:\Omega_2\rightarrow\Omega$ be the exponential map of $g$. That is, for all $(x,y)$, the curve $t\mapsto\opExp(x,ty)$ is the unique geodesic in $\Omega$ leaving the point $x$ in the direction of the vector $y$.
\noskipproclaim{Proposition \nextprocno}
$$
[\opExp] \in \overline{\Cal{R}}^{-1}[[X,Y]],\eqnum{\nexteqnno[TaylorSeriesOfExp]}
$$
and
$$
[\opExp] = X + Y + O(\|X,Y\|^3).\eqnum{\nexteqnno[FirstTermsOfTaylorSeriesOfExp]}
$$
\endproclaim
\proclabel{TaylorSeriesOfExp}
\proof For any function $\phi$ of $s$ and $t$, and for all $k$, let $[\phi]_{\infty,k}$ denote its Taylor series up to order $k$ in $t$. Likewise, for any formal series $\Phi$ in $S$ and $T$, let $[\Phi]_{\infty,k}$ denote its partial sum up to order $k$ in $T$. Now define $E(x,y,s,t):=\opExp(sx,ty)$. By definition, for every point $(x,y)\in\Omega_2$ and for all $s$,
$$\triplealign{
E(x,y,s,0) &= \opExp(sx,0) &=sx,\ \&\cr
\partial_t E(x,y,s,0) &= \partial_t \opExp(sx,ty)|_{t=0}&=y.\cr}
$$
so that $[E]_{\infty,1}=SX+TY$, which belongs to $\overline{\Cal{R}}[X,Y][[S,T]]$. We now claim that the partial sum $[E]_{\infty,k}$ belongs to $\overline{\Cal{R}}[X,Y][[S,T]]$ for all $k$. Indeed, suppose that this holds for $k$. Observe that
$$\eqalign{
[\Gamma(E)(\partial_t E,\partial_t E)]_{\infty,k-1}
&= [\Gamma([E]_{\infty,k-1})(\partial_T[E]_{\infty,k},\partial_T[E]_{\infty,k})]_{\infty,k-1},\cr}
$$
where $\partial_T$ denotes formal partial differentiation with respect to $T$. Since $[\Gamma]$ belongs to $\overline{\Cal{R}}[X]$, it follows by the inductive hypothesis that $[\Gamma(E)(\partial_t E,\partial_t E)]_{\infty,k-1}$ belongs to $\overline{\Cal{R}}[X,Y][[S,T]]$. However, by the geodesic equation,
$$
\partial_T^2[E]_{\infty,k+1} = [\partial_t^2 E]_{\infty,k-1} = -[\Gamma(E)(\partial_t E,\partial_t E)]_{\infty,k-1} \in \overline{\Cal{R}}[X,Y][[S,T]],
$$
and the claim now follows by induction. In particular $[E]$ belongs to $\overline{\Cal{R}}[X,Y][[S,T]]$ and the first assertion follows by Proposition \procref{TestForBeingCurvatureSeries}. Finally, since $[\opExp]-(X+Y)\in\Cal{R}^{-1}[[X,Y]]$, its lowest degree term has degree at least $3$ in $X$ and $Y$, thus proving the second assertion. This completes the proof.\qed
\medskip
Let $\opTr:\Omega_2\times\Bbb{R}^{n+1}\rightarrow\opEnd(\Bbb{R}^{n+1})$ be such that for all $(x,y)$ and for every vector $U$, $\opTr(x,y)U$ is the parallel transport of $U$ from the point $x$ to the point $\opExp(x,y)$ along the geodesic $t\mapsto\opExp(x,ty)$.
\noskipproclaim{Proposition \nextprocno}
$$
[\opTr] \in\overline{\Cal{R}}^0[[X,Y]],\eqnum{\nexteqnno[TaylorSeriesOfT]}
$$
and
$$
[\opTr] = I + O(\|X,Y\|^2).\eqnum{\nexteqnno[FirstTermsOfTaylorSeriesOfT]}
$$
\endproclaim
\proclabel{TaylorSeriesOfParallelTransport}
\proof As before, for any function $\phi$ of $s$ and $t$, and for all $k$, let $[\phi]_{\infty,k}$ denote its Taylor series up to order $k$ in $t$. Likewise, for any formal series $\Phi$ in $S$ and $T$, let $[\Phi]_{\infty,k}$ denote its partial sum up to order $k$ in $T$. Define $E(x,y,s,t):=\opExp(sx,ty)$ and $F(x,y,s,t):=\opTr(sx,ty)$. By definition, for every point $(x,y)\in\Omega_2$ and for all $s$,
$$
F(x,y,s,0) = \opTr(sx,0) = I,
$$
so that $[F]_{\infty,0}=I$, which belongs to $\overline{\Cal{R}}[X,Y][[S,T]]$. We now claim that the partial sum $[F]_{\infty,k}$ belongs to $\overline{\Cal{R}}[X,Y][[S,T]]$ for all $k$. Indeed, suppose that this holds for $k$. Observe that
$$\eqalign{
[\Gamma(E)(\partial_tE,F)]_{\infty,k}
&= [\Gamma([E]_{\infty,k})(\partial_T[E]_{\infty,k+1},[F]_{\infty,k})]_{\infty,k},\cr}
$$
where $\partial_T$ denotes formal partial differentiation with respect to $T$. Since $[\Gamma]$ belongs to $\overline{\Cal{R}}[X]$ and since $[\opExp]$ belongs to $\overline{\Cal{R}}[X,Y]$, it follows by the inductive hypothesis that $[\Gamma(E)(\partial_tE,F)]_{\infty,k}$ also belongs to $\overline{\Cal{R}}[X,Y][[S,T]]$. However, by the parallel transport equation
$$
\partial_T[F]_{\infty,k+1} = [\partial_tF]_{\infty,k} = -[\Gamma(E)(\partial_t E,F)]_{\infty,k} \in \overline{\Cal{R}}[X,Y][[S,T]],
$$
and the claim now follows by induction. In particular, $F$ belongs to $\overline{\Cal{R}}[X,Y][[S,T]]$ and the first assertion follows by Proposition \procref{TestForBeingCurvatureSeries}. Finally, since $[F]-I\in\Cal{R}^0[[X,Y]]$, its lowest degree term has degree at least $2$ in $X$ and $Y$, thus proving the second assertion. This completes the proof.\qed
\medskip
Finally, we consider higher order iterates of the exponential map and the parallel transport. Thus, for all $n$, define $\Omega_{n+1}\subseteq(\Bbb{R}^{m+1})^{n+1}$ by,
$$
\Omega_{n+1} = \left\{ (x_1,...,x_{n+1})\ |\ \|x_1\| + ... + \|x_{n+1}\| < 1\right\},
\eqnum{\nexteqnno[AdvancedDomainOfExponentialMap]}
$$
and define the sequences of functions $(\opExp_n)$ and $(\opTr_n)$ such that
$$\eqalign{
\opExp_1(x_1,x_2) &:=\opExp(x_1,x_2),\ \&\cr
\opTr_0(x_1) &:=\opId,\cr}
\eqnum{\nexteqnno[InductiveDefinitionOfExpAndTStepI]}
$$
and, for all $n$,
$$\eqalign{
\opExp_n(x_1,...,x_{n+1}) \hfill&:= \opExp(\opExp_{n-1}(x_1,...,x_n),\opTr_{n-1}(x_1,...,x_n)x_{n+1}),\ \&\hfill\cr
\opTr_n(x_1,...,x_{n+1})U \hfill&:= \opTr(\opExp_{n-1}(x_1,...,x_n),\opTr_{n-1}(x_1,...,x_{n-1})U).\hfill\cr}
\eqnum{\nexteqnno[InductiveDefinitionOfExpAndTStepII]}
$$
\proclaim{Proposition \nextprocno}
\noindent For all $n$,
$$\eqalign{
[\opExp_n] &\in \overline{\Cal{R}}^{-1}[[X_1,...,X_{n+1}]],\ \&\cr
[\opTr_n] &\in \overline{\Cal{R}}^{0}[[X_1,...,X_{n+1}]],\cr}
\eqnum{\nexteqnno[TaylorSeriesOfExtendedFunctions]}
$$
and
$$\eqalign{
[\opExp_n] &= X_1 + ... + X_{n+1} + O(\|X_1,...,X_{n+1}\|^3),\ \&\cr
[\opTr_n] &= I + O(\|X_1,...,X_{n+1}\|^2).\cr}
\eqnum{\nexteqnno[FirstTermsOfTaylorSeriesOfExtendedFunctions]}
$$
\endproclaim
\proclabel{TaylorSeriesOfExtendedFunctions}
\proof This follows by induction using Propositions \procref{TaylorSeriesOfExp} and \procref{TaylorSeriesOfParallelTransport} and the recursive definitions of $(\opExp_n)$ and $(\opTr_n)$.\qed
\newhead{Taylor Series of Functions Derived From Immersions}[TaylorSeriesOfFunctionsDerivedFromImmersions]
\newsubhead{Graphs Over Spheres}[GraphsOverSpheres]
Let $S^m$ be the unit sphere in $\Bbb{R}^{m+1}$ and let $\overline{\nabla}$, $\overline{\opHess}$ and $\overline{\Delta}$ be respectively its gradient, Hessian and Laplace operators with respect to the standard euclidean metric. For $t\in]0,\infty[$, which we think of as a scale parameter, and for $f\in C^0(S^m)$, consider the function $e(t,f):S^m\rightarrow\Bbb{R}^{m+1}$ given by
$$
e(t,f)(x) := t(1+t^2 f(x)) x.\eqnum{\nexteqnno[GraphOverSphere]}
$$
Heuristically, $e(t,f)$ is an immersed sphere of radius approximately $t$ centred on the origin. For all $k$, let $J:=J^kS^m$ denote the bundle of $k$-jets over $S^m$, and for a function $f\in C^k(S^m)$ and a point $x\in S^m$, denote by $f_x$ its $k$-jet at $x$, where the order $k$ of the jet should hopefully be clear from the context. Define the functions $N:]0,\infty[\times J^1S^m\rightarrow S^m$ and $H:]0,\infty[\times J^2S^m\rightarrow\Bbb{R}$ such that for all $t\in]0,\infty[$ and for all $f_x\in JS^m$, $N(t,f_x)$ and $H(t,f_x)$ are respectively the  outward-pointing unit normal of the immersion $e(t,f)$ at the point $e(t,f)(x)$ and its mean curvature at that point, both with respect to the metric $g$. It is important to note that both $N$ and $H$ are actually smooth functions defined over finite-dimensional domains and may both be expressed explicitly in terms of (rather complicated) formulae involving $g$. We have chosen to define these functions in this apparently roundabout manner in order to emphasise their clear geometric meanings.
\par
We are interested in the Taylor series of $N(t,f_x)$ and $H(t,f_x)$ in $t$ about $0$. To this end, we first introduce the following auxiliary functions. Define $r:\Bbb{R}^{m+1}\rightarrow[0,\infty[$ and $x:\Bbb{R}^{m+1}\setminus\left\{0\right\}\rightarrow S^m$ by
$$\eqalign{
r(y)&:=\|y\|, \&\cr
x(y)&:=y/r.\cr}\eqnum{\nexteqnno[DefinitionsOfXAndR]}
$$
Given $f\in C^1(S^m)$, define $\hat{f}:]0,\infty[\times(\Bbb{R}^{m+1}\setminus\left\{0\right\})\rightarrow\Bbb{R}$ by
$$
\hat{f}(t,y) := r - t(1 + t^2 f(x)).\eqnum{\nexteqnno[FormulaForF]}
$$
Observe that the image of $e(t,f)$ coincides with the level set of $\hat{f}$ at height $0$. Furthermore, for every point $y$ in this level set, $\nabla\hat{f}(y)$ is orthogonal to this level set with respect to the metric $g$.
\noskipproclaim{Proposition \nextprocno}
$$
\nabla \hat{f}(t,y) = \frac{y}{r} - \frac{t}{r}B(y)t^2\overline{\nabla} f(x).\eqnum{\nexteqnno[GradientOfF]}
$$
\endproclaim
\proclabel{GradientOfF}
\proof The gradient of $\hat{f}$ with respect to the euclidean metric is
$$
D\hat{f}(t,y) = \frac{y}{r} - \frac{t^3}{r}\overline{\nabla}f(x).
$$
However, for all vectors $U$ in $\Bbb{R}^{m+1}$,
$$
d\hat{f}(t,y)(U) = \langle D\hat{f}(t,y),U\rangle = \langle A(y)B(y)D\hat{f}(t,y),U\rangle = g(B(y)D\hat{f}(t,y),U),
$$
so that the gradient of $\hat{f}$ with respect to $g$ is $\nabla\hat{f}(t,y) = B(y)D\hat{f}(t,y)$, and since $B(y)y=y$ for all $y$, the result follows.\qed
\medskip
We now invert the situation and consider both $r$ and $y$ as functions of $t$ and $x$, so that
$$\eqalign{
r(t,x) &:= t(1 + t^2f(x)),\ \&\cr
y(t,x) &:= t(1 + t^2f(x))x.\cr}\eqnum{\nexteqnno[FormulaeForRAndY]}
$$
We define
$$
\hat{N}(t,f_x) := \nabla\hat{f}(t,y) = \frac{y}{r} - \frac{t}{r}B(y)t^2\overline{\nabla}f(x),\eqnum{\nexteqnno[FormulaForNHat]}
$$
so that we obtain the following formula for $N$.
$$
N(t,f_x) := \frac{1}{\|\hat{N}(t,x)\|_g}\hat{N}(t,f_x),\eqnum{\nexteqnno[FormulaForN]}
$$
where $\|\cdot\|_g$ here denotes the norm with respect to the metric $g$.
\par
It will also be necessary to extend $e$, $N$ and $H$ to allow for variations of the centre of the immersed sphere. Thus, for $t\in]0,\infty[$, for $y\in\Bbb{R}^{m+1}$, and for $f\in C^0(S^m)$, define $e(t,y,f):S^m\rightarrow\Bbb{R}^{m+1}$ by
$$
e(t,y,f)(x) := \opExp(ty,t(1+t^2f(x))x),\eqnum{\nexteqnno[FormulaForImmersionVariableCentre]}
$$
so that, heuristically, $e(t,y,f)$ is an immersed sphere of radius approximately $t$ with centre displaced to the point $y$. Define $N:]0,\infty[\times\Bbb{R}^{m+1}\times J^1S^m\rightarrow S^m$ and $H:]0,\infty[\times\Bbb{R}^{m+1}\times J^2S^m\rightarrow\Bbb{R}$ as before. Observe, in particular, that $e(t,0,f)=e(t,f)$, $N(t,0,f_x)=N(t,f_x)$ and $H(t,0,f_x)=H(t,f_x)$.
\newsubhead{The Taylor Series of the Unit Normal Vector}[TheTaylorSeriesOfTheUnitNormalVector]
We now study the Taylor series of the scale-dependent functions introduced in Section \subheadref{GraphsOverSpheres}. In particular, we are interested in how the different terms in these series contribute to the exponent of $T$. To this end, we extend the formalism developed in Sections \subheadref{CurvatureTensors} and \subheadref{CurvaturePolynomials}. Thus, for a vector $\underline{X}:=(X_1,...,X_n)$ of formal variables taking values in $\Bbb{R}^{m+1}$, consider the set of formal polynomials
$$
\left\{X^a_i\delta_{ab}P^b_{j_1...j_k}(\underline{X})\ |\ P\in\overline{\Cal{R}}[\underline{X}]\right\},
\eqnum{\nexteqnno[GeneratorsOfCurvaturePolynomials]}
$$
where $X^i_j$ denotes the $i$'th component of the vector $X_j$. Let $\Cal{Q}[\underline{X}]$ be the vector space with basis the set of all tensor products of elements of this set. We call $\Cal{Q}[\underline{X}]$ the space of {\bf curvature polynomials} of the second kind. For all $k\in\Bbb{N}$ and for all $\underline{r}:=(r_1,...,r_n)\in\Bbb{N}^n$, denote by $\Cal{Q}^k_{\underline{r}}[\underline{X}]$ the subspace consisting of those elements which are contravariant of order $k$ and which are homogeneous of degree $r_i$ in the variable $X_i$ for all $i$. When $Q\in\Cal{R}^k_{\underline{r}}[\underline{X}]$, we say that it has {\bf order} $k$ and {\bf degree} $\underline{r}$. Finally, denote
$$
\overline{\Cal{Q}}[\underline{X}] := \Cal{Q}[\underline{X}]\oplus\langle 1\rangle.\eqnum{\nexteqnno[ExtendedCurvaturePolynomialsOfTheSecondKind]}
$$
We also call elements of $\overline{\Cal{Q}}[\underline{X}]$ {\bf curvature polynomials} of the second kind.
\par
Now let $A$ be an algebra graded by $\Bbb{N}^k$ for some $k$. Let $A[T]$ be the algebra of polynomials over $\Bbb{R}$ with coefficients in $A$. For a given weight $\underline{w}:=(w_1,...,w_k)\in\Bbb{N}^k$, let $A[T]_{\underline{w}}$ be the subalgebra of $A[T]$ consisting of those polynomials whose coefficients of degree $m$ are elements of $\oplus_{\langle\underline{w},\underline{i}\rangle=m}A_{\underline{i}}$ for all $m$. Likewise, let $A[[T]]$ be the algebra of formal power series over $\Bbb{R}$ with coefficients in $A$, and for $\underline{w}\in\Bbb{N}^k$, let $A[[T]]_{\underline{w}}$ be the subalgebra of $A[[T]]$ consisting of those formal power series all of whose partial sums are elements of $A[T]_{\underline{w}}$.
\par
Now let $\Bbb{R}[F]$ be the algebra of formal polynomials in the variable $F$. Consider a smooth function $\phi:[0,\infty[\times J^kS^m\rightarrow\Bbb{R}$ which only depends on the metric $g$ and the jet $f_x$. For such a function, the statement that $[\phi]$ belongs to $\Bbb{R}_*[F]\otimes\overline{\Cal{Q}}_{*,*}[X,\overline{\nabla}F][[T]]_{(2,1,2)}$, for example, means that its Taylor series in $t$ about $0$ takes the form
$$
\phi(t,f_x) \sim \sum_{m=0}^\infty t^m \sum_{\langle\underline{i},(2,1,2)\rangle=m}
\sum_{\alpha}P_{\underline{i},\alpha}(f(x))Q_{\underline{i},\alpha}(x,\overline{\nabla} f(x)),
\eqnum{\nexteqnno[MeaningOfFormalism]}
$$
where, for all $\underline{i}:=(i_1,i_2,i_3)$ and for all $\alpha$, $P_{\underline{i},\alpha}$ is a polynomial of degree $i_1$ and $Q_{\underline{i},\alpha}$ is a curvature polynomial of the second kind of order $0$ and degree $(i_2,i_3)$. We leave the reader to interpret the meanings of other tensor products of spaces of formal polynomials. Importantly, this notation emphasises that all terms in $X$ carry weight $1$ in $T$ whilst all terms in $F$ and $\overline{\nabla}F$ carry weight $2$. This behaviour will be common to all series studied in the sequel.
\proclaim{Proposition \nextprocno}
\noindent For all real $\alpha$,
$$
[(r/t)^\alpha]\in\Bbb{R}_*[F][[T]]_2.\eqnum{\nexteqnno[TaylorSeriesOfROverT]}
$$
\endproclaim
\proclabel{TaylorSeriesOfROverT}
\proof By definition, $[(r/t)]=[1+t^2f]$ belongs to $\Bbb{R}_*[F][[T]]_2$. Since $[(r/t)](0)= (r/t)(0)=1$, the result follows by Proposition \procref{PowersOfCurvatureSeriesAreAlsoCurvatureSeries}.\qed
\proclaim{Proposition \nextprocno}
\noindent For all real $\alpha$,
$$
[\|\hat{N}(t,f_x)\|^\alpha]\in\Bbb{R}_*[F]\otimes\overline{\Cal{Q}}^0_{*,*}[X,\overline{\nabla}F][[T]]_{(2,1,2)},
\eqnum{\nexteqnno[TaylorSeriesOfNormN]}
$$
and,
$$
[\|\hat{N}(t,f_x)\|^\alpha]=1+O(T^4).
\eqnum{\nexteqnno[ShortTaylorSeriesOfNormN]}
$$
\endproclaim
\proclabel{TaylorSeriesOfNormN}
\proof Using \eqnref{MetricIsExponential}, \eqnref{EqnMetric} and \eqnref{FormulaForNHat}, we obtain, for all $t$ and for all $x$,
$$
\|\hat{N}(t,f_x)\|_g^2 = 1 + \left(\frac{t}{r}\right)^2\langle B(y)t^2\overline{\nabla} f, t^2\overline{\nabla} f\rangle.
$$
However, by Propositions \procref{TaylorSeriesOfAAndB} and \procref{TaylorSeriesOfROverT}, $[B(y)]=[B((r/t)(tx))]$ belongs to $\Bbb{R}_*[F]\otimes\overline{\Cal{R}}_*[X][[T]]_{(2,1)}$, so that $[\langle B(y)t^2\overline{\nabla}f,t^2\overline{\nabla}f\rangle]$ belongs to $\Bbb{R}_*[F]\otimes\overline{\Cal{Q}}_{*,*}[X,\overline{\nabla} F][[T]]_{(2,1,2)}$. It follows by Proposition \procref{TaylorSeriesOfROverT} again that $\|\hat{N}(t,f_x)\|_g^2$ belongs to $\Bbb{R}_*[F]\otimes\overline{\Cal{Q}}_{*,*}[X,\overline{\nabla}F][[T]]_{(2,1,2)}$, and, since the first term in this series equals $1$, the first assertion follows by Proposition \procref{PowersOfCurvatureSeriesAreAlsoCurvatureSeries}. Finally, since $[(t/r)^2]$ has order $0$ in $T$ and since $[\langle B(y) t^2\overline{\nabla} f,t^2\overline{\nabla}f\rangle]$ has order $4$ in $T$, we see that $\|\hat{N}(t,f_x)\|^2_g=1+O(T^4)$, and the second assertion follows by Proposition \procref{PowersOfCurvatureSeriesAreAlsoCurvatureSeries} again. This completes the proof.\qed
\noskipproclaim{Proposition \nextprocno}
$$
N(t,f_x) = \Phi_1(t,f_x)x + \Phi_2(t,f_x),\eqnum{\nexteqnno[TaylorSeriesOfNPartI]}
$$
where
$$\eqalign{
[\Phi_1] &\in \Bbb{R}_*[F]\otimes\overline{\Cal{Q}}^0_{*,*}[X,\overline{\nabla} F][[T]]_{(2,1,2)},\ \&\cr
[\Phi_2] &\in \Bbb{R}_*[F]\otimes\overline{\Cal{Q}}^0_{*,*}[X,\overline{\nabla}F]\otimes
\overline{\Cal{R}}^{-1}_{*,*}[X,\overline{\nabla}F][[T]]_{(2,1,2,1,2)}.\cr}
\eqnum{\nexteqnno[TaylorSeriesOfNPartII]}
$$
Furthermore
$$
[N(t,f_x)] = X - T^2\overline{\nabla}F + O(T^4).\eqnum{\nexteqnno[FirstTermsInTaylorSeriesOfN]}
$$
\endproclaim
\proclabel{TaylorSeriesOfN}
\proof As in the proof of Proposition \procref{TaylorSeriesOfNormN}, $[B(y)]$ belongs to $\Bbb{R}_*[F]\otimes\overline{\Cal{R}}_*[X][[T]]_{(2,1)}$ and so $[B(y)t^2\overline{\nabla}f]$ belongs to $\Bbb{R}_*[F]\otimes\overline{\Cal{R}}_{*,*}[X,\overline{\nabla}F][[T]]_{(2,1,2)}$. By Proposition \procref{TaylorSeriesOfROverT}, the series $[(t/r)B(y)t^2\overline{\nabla}f]$ also belongs to $\Bbb{R}_*[F]\otimes\overline{\Cal{R}}_{*,*}[X,\overline{\nabla}F][[T]]_{(2,1,2)}$, and the result now follows by \eqnref{FormulaForNHat}, \eqnref{FormulaForN} and Proposition \procref{TaylorSeriesOfNormN}.\qed
\noskipproclaim{Proposition \nextprocno}
$$
N(t,y,f_x) = \Phi_1(t,y,f_x)x + \Phi_2(t,y,f_x),\eqnum{\nexteqnno[TaylorSeriesOfNVariableCentrePartI]}
$$
where
$$\eqalign{
[\Phi_1] &\in \Bbb{R}_*[F]\otimes\overline{\Cal{Q}}^0_{*,*,*}[X,Y,\overline{\nabla}F][[T]]_{(2,1,1,2)},\ \&\cr
[\Phi_2] &\in \Bbb{R}_*[F]\otimes\overline{\Cal{Q}}^0_{*,*,*}[X,Y,\overline{\nabla}F]\otimes
\overline{\Cal{R}}^{-1}_{*,*,*}[X,Y,\overline{\nabla}F][[T]]_{(2,1,1,2,1,1,2)}.\cr}
\eqnum{\nexteqnno[TaylorSeriesOfNVariableCentrePartII]}
$$
Furthermore
$$
[N(t,f_x)] = X + O(T^2).\eqnum{\nexteqnno[FirstTermsInTaylorSeriesOfNVariableCentre]}
$$
\endproclaim
\proclabel{TaylorSeriesOfNVariableCentre}
\proof This Taylor series is obtained from Proposition \procref{TaylorSeriesOfN} by substituting for every generator $R_{i_1i_2i_3}{}^j{}_{;i_4...i_{k+3}}$ of $\Cal{R}$, its own Taylor series in $t$ about $0$,
$$
[R_{i_1i_2i_3}{}^j{}_{;i_4...i_{k+3}}]=\sum_{m=0}^\infty\frac{1}{m!}T^m R_{i_1i_2i_3}{}^j{}_{;i_4...i_{k+3+m}}Y^{i_{k+3+1}}...Y^{i_{k+3+m}}.
$$
The result follows.\qed
\newsubhead{Normal Variation of Spheres}[NormalVariation]
We extend $e$ further in order to study variations of the base point, of the displacement of the centre, and of the immersion itself. Thus, for $t\in]0,\infty[$, for $y,z,w\in\Bbb{R}^{m+1}$ and for $f,g\in C^0(S^m)$, consider the function $e(t,y,z,w,f,g):S^m\rightarrow\Bbb{R}^{m+1}$ given by
$$
e(t,y,z,w,f,g)(x) := \opExp_2(z,t(y+w),t(1+t^2(f(x)+g(x)))x),\eqnum{\nexteqnno[MostGeneralVariation]}
$$
and define $P,Q:]0,\infty[\times\Bbb{R}^{m+1}\times J^0S^m\rightarrow\opEnd(\Bbb{R}^{m+1})$ and $R:]0,\infty[\times\Bbb{R}^{m+1}\times J^0S^m\rightarrow\Bbb{R}^{m+1}$ by
$$\eqalign{
P(t,y,f_x) &:= \partial_z e(t,y,0,0,f,0)(x),\cr
Q(t,y,f_x) &:= \partial_w e(t,y,0,0,f,0)(x),\ \&\cr
R(t,y,f_x) &:= \partial_g e(t,y,0,0,f,0)(x).\cr}
\eqnum{\nexteqnno[DefnOfPQR]}
$$
Heuristically, for any given vectors, $U$ and $V$, and for any given function, $g$, the vectors $P(t,y,f_x)U$, $Q(t,y,f_x)V$ and $R(t,y,f_x)g_x$ measure the respective infinitesimal variations of the immersion $e(t,y,f)$ at the point $e(t,y,f)(x)$ arising from infinitesimal perturbations of the base point, of the displacement of the centre, and of the immersion itself in the directions of $U$, $tV$ and $t^3g$ respectively.
\par
Now define $p,q:]0,\infty[\times\Bbb{R}^{m+1}\times J^0 S^m\rightarrow\Bbb{R}^{m+1}$ and $r:]0,\infty[\times\Bbb{R}^{m+1}\times J^0S^m\rightarrow\Bbb{R}$ by
$$\eqalign{
\langle p(t,y,f_x),U\rangle &:= \langle A(e(t,y,f_x))P(t,y,f_x)U,N(t,y,f_x)\rangle,\cr
\langle q(t,y,f_x),V\rangle &:= \langle A(e(t,y,f_x))Q(t,y,f_x)V,N(t,y,f_x)\rangle,\ \&\cr
r(t,y,f_x)g &:= \langle A(e(t,y,f_x))R(t,y,f_x)g,N(t,y,f_x)\rangle.\cr}
\eqnum{\nexteqnno[DefnOfpqr]}
$$
Heuristically, $p$, $q$ and $r$ measure the normal components of the above infinitesimal variations.
\noskipproclaim{Proposition \nextprocno}
$$
p(t,y,f_x) = \Phi_1(t,y,f_x)x + \Phi_2(t,y,f_x),\eqnum{\nexteqnno[TaylorSeriesOfLittlePPartI]}
$$
where
$$\eqalign{
[\Phi_1] &\in \Bbb{R}_*[F]\otimes\overline{\Cal{Q}}^0_{*,*,*}[X,Y,\overline{\nabla}F]
\otimes\overline{\Cal{R}}^0_{*,*,*}[X,Y,\overline{\nabla}F][[T]]_{(2,1,1,2,1,1,2)}, \&\cr
[\Phi_2] &\in \Bbb{R}_*[F]\otimes\overline{\Cal{Q}}^0_{*,*,*}[X,Y,\overline{\nabla}F]
\otimes\overline{\Cal{R}}^{-1}_{*,*,*}[X,Y,\overline{\nabla}F][[T]]_{(2,1,1,2,1,1,2)}.\cr}
\eqnum{\nexteqnno[TaylorSeriesOfLittleP]}
$$
Furthermore
$$
[p] = \langle X, \cdot\rangle + O(T^2).\eqnum{\nexteqnno[FirstTermsInTaylorSeriesOfLittleP]}
$$
\endproclaim
\proclabel{TaylorSeriesOfLittleP}
\proof Let $\partial_Z$ denote the formal partial derivative with respect to the variable $Z$. In particular, $[\partial_z\opExp_2(z,y,x)|_{z=0}]=\partial_Z[\opExp_2(z,y,x)]|_{Z=0}$. However, by Proposition \procref{TaylorSeriesOfExtendedFunctions},
$$\eqalign{
\partial_Z[\opExp_2(z,y,x)]|_{Z=0} &\in \overline{\Cal{R}}[X,Y],\ \&\cr
\partial_Z[\opExp_2(z,y,x)]|_{Z=0} &= I + O(\|X,Y\|^2).}
$$
Substituting $ty$ and $t(1+t^2f(x))x$ for $y$ and $x$ respectively therefore yields
$$\eqalign{
[\partial_ze(z,ty,t(1+t^2f(x))x)|_{z=0}] &\in \Bbb{R}_*[F]\otimes\overline{\Cal{R}}_{*,*}[X,Y][[T]]_{(2,1,1)},\ \&\cr
[\partial_ze(z,ty,t(1+t^2f(x))x)|_{z=0}] &= I+O(T^2),\cr}
$$
so that
$$\eqalign{
[P]&\in\Bbb{R}_*[F]\otimes\overline{\Cal{R}}^0_{*,*}[X,Y][[T]]_{(2,1,1)},\ \&\cr
[P]&=I+O(T^2).\cr}
$$
The result now follows from the self-adjointness of $\Cal{R}$ (Proposition \procref{TheSpaceOfCurvaturePolynomialsIsSelfAdjoint}) and Propositions \procref{TaylorSeriesOfAAndB} and \procref{TaylorSeriesOfNVariableCentre}.\qed
\noskipproclaim{Proposition \nextprocno}
$$
q(t,y,f_x) = t\Phi_1(t,y,f_x)x + t\Phi_2(t,y,f_x),\eqnum{\nexteqnno[TaylorSeriesOfLittleQPartI]}
$$
where
$$\eqalign{
[\Phi_1] &\in \Bbb{R}_*[F]\otimes\overline{\Cal{Q}}^0_{*,*,*}[X,Y,\overline{\nabla}F]
\otimes\overline{\Cal{R}}^0_{*,*,*}[X,Y,\overline{\nabla}F][[T]]_{(2,1,1,2,1,1,2)},\ \&\cr
[\Phi_2] &\in \Bbb{R}_*[F]\otimes\overline{\Cal{Q}}^0_{*,*,*}[X,Y,\overline{\nabla}F]
\otimes\overline{\Cal{R}}^{-1}_{*,*,*}[X,Y,\overline{\nabla}F][[T]]_{(2,1,1,2,1,1,2)}.\cr}
\eqnum{\nexteqnno[TaylorSeriesOfLittleQ]}
$$
Furthermore
$$
[q] = T\langle X, \cdot\rangle + O(T^3).\eqnum{\nexteqnno[FirstTermsInTaylorSeriesOfLittleQ]}
$$
\endproclaim
\proclabel{TaylorSeriesOfLittleQ}
\proof Let $\partial_W$ denote the formal partial derivative with respect to the variable $W$. In particular, $[\partial_w\opExp(y+tw,x)|_{w=0}]=\partial_W[\opExp(y+tw,x)]|_{W=0}$. However, by Proposition \procref{TaylorSeriesOfExtendedFunctions},
$$\eqalign{
\partial_W[\opExp(y+tw,x)]|_{W=0} &\in T\overline{\Cal{R}}[X,Y],\ \&\cr
\partial_W[\opExp(y+tw,x)]|_{W=0} &= TI + TO(\|X,Y\|^2).}
$$
Substituting $ty$ and $t(1+t^2f(x))x$ for $y$ and $x$ respectively therefore yields
$$\eqalign{
[\partial_we(z,ty,t(1+t^2f(x))x)|_{w=0}] &\in T\Bbb{R}_*[F]\otimes\overline{\Cal{R}}_{*,*}[X,Y][[T]]_{(2,1,1)},\ \& \cr
[\partial_we(z,ty,t(1+t^2f(x))x)|_{w=0}] &= TI+O(T^3),\cr}
$$
so that
$$\eqalign{
[Q]&\in T\Bbb{R}_*[F]\otimes\overline{\Cal{R}}^0_{*,*}[X,Y][[T]]_{(2,1,1)},\cr
[Q]&=TI+O(T^3).\cr}
$$
The result now follows from the self-adjointness of $\Cal{R}$ (Proposition \procref{TheSpaceOfCurvaturePolynomialsIsSelfAdjoint}) and Propositions \procref{TaylorSeriesOfAAndB} and \procref{TaylorSeriesOfNVariableCentre}.\qed
\noskipproclaim{Proposition \nextprocno}
$$
[r] \in T^3\Bbb{R}_*[F]\otimes\overline{\Cal{Q}}^0_{*,*,*}[X,Y,\overline{\nabla}F][[T]]_{(2,1,1,2)},\eqnum{\nexteqnno[TaylorSeriesOfLittleR]}
$$
and
$$
[r] = T^3 + O(T^7).\eqnum{\nexteqnno[FirstTermsInTaylorSeriesOfLittleR]}
$$
\endproclaim
\proclabel{TaylorSeriesOfLittleR}
\proof Consider first the case where $Y=0$ and observe that
$$
e(t,0,0,0,f,g) = \opExp(t(1+t^2(f(x)+g(x)))x).
$$
In particular, since $\Omega$ is an exponential chart,
$$
\partial_g\opExp(t(1+t^2(f(x)+g(x)))x)|_{g=0} = t^3x,
$$
so that $R(t,0,f_x)=t^3x$. Thus, by \eqnref{FormulaForNHat} and \eqnref{FormulaForN}, since $A(y)x=x$ and since $\langle x,\overline{\nabla}f(x)\rangle=0$,
$$
r(t,0,f_x) = t^3 \langle x,N(t,0,f_x)\rangle = t^3\|\hat{N}(t,x)\|_g^{-1}.
$$
The result for $Y=0$ now follows by Proposition \procref{TaylorSeriesOfNormN}. The result for general $Y$ follows by substituting for every generator $R_{i_1i_2i_3}{}^j_{;i_4...i_{k+3}}$ of $\Cal{R}$ its own Taylor series in $Y$ about $0$, as in the proof of Proposition \procref{TaylorSeriesOfNVariableCentre}.\qed
\newsubhead{The Taylor Series of the Mean Curvature}[TheTaylorSeriesOfTheMeanCurvature]
We end this section by determining the Taylor series of the mean curvature function. First recall that (c.f. \cite{SmiCCH}),
$$\eqalign{
H(t,Y,f_x) &\sim \frac{1}{t}\left(1 - \frac{t^2}{3}\opRic_{pq}x^px^q - \frac{t^2}{n}(n+\overline{\Delta})f - \frac{t^3}{4}\opRic_{pq;r}x^px^qx^r\right.\cr
&\qquad -\left.\frac{t^3}{3}\opRic_{pq;r}x^px^pY^r - \frac{t^4}{4}\opRic_{pq;rs}x^px^qx^rY^s + t^4F(f_x) + O(t^5)\right),\cr}
\eqnum{\nexteqnno[FirstAsymptoticFormulaForDisplacedCurvature]}
$$
where $F$ is a curvature polynomial. More generally,
\noskipproclaim{Proposition \nextprocno}
$$
H(t,y,f_x) = \frac{1}{t}\opTr(\Phi_1) + \opTr(\Phi_2) + \frac{1}{t}\opTr(\Phi_3 t^2\overline{\opHess}(f)\circ\pi) +
\frac{1}{t}\langle\Phi_5,t^2\overline{\opHess}(f)\circ\pi\rangle,\eqnum{\nexteqnno[TaylorSeriesOfHPartI]}
$$
where
$$\eqalign{
\Phi_1,\Phi_2,\Phi_3 &\in \Bbb{R}_*[F]\otimes\overline{\Cal{Q}}^0_{*,*}[X,Y,\overline{\nabla} F]\otimes\overline{\Cal{R}}^0_{*,*}[X,Y,\overline{\nabla}F][[T]]_{(2,1,1,2,1,1,2)},\ \&\cr
\Phi_4 &\in \Bbb{R}_*[F]\otimes\overline{\Cal{Q}}^2_{*,*}[X,Y,\overline{\nabla} F][[T]]_{(2,1,1,2)}\cr}
\eqnum{\nexteqnno[TaylorSeriesOfHPartII]}
$$
\endproclaim
\proclabel{TaylorSeriesOfH}
\proof We first consider the case where $Y=0$. Recall that
$$
H = \frac{1}{\|\nabla\hat{f}\|_g}\Delta\hat{f} - \frac{1}{\|\nabla\hat{f}\|_g^3} g(\nabla_{\nabla\hat{f}}\nabla \hat{f},\nabla\hat{f}),\eqnum{\nexteqnno[FormulaForMeanCurvature]}
$$
where $\Delta$ here denotes the Laplace operator of the metric $g$. Furthermore, by \eqnref{GradientOfF},
$$
\nabla\hat{f} = \frac{1}{r}\left(y - tB(y)t^2\overline{\nabla} f(x)\right).\eqnum{\nexteqnno[GradientOfFAgain]}
$$
Now observe that, for all vectors $U$,
$$
D_U\overline{\nabla}f(x) = \frac{1}{r}\overline{\opHess}(f)\circ\pi(U) + \frac{1}{r}\langle U,\overline{\nabla}f(x)\rangle\frac{y}{r},
$$
where $\pi$ is the orthogonal projection along $x$. Differentiating \eqnref{GradientOfFAgain}, therefore yields, for all $U$,
$$\eqalign{
\nabla_U\nabla\hat{f}
&=\frac{1}{t}\left(\frac{t}{r}\right)U
-\frac{1}{t}\left(\frac{t}{r}\right)\left\langle U,\frac{y}{r}\right\rangle\frac{y}{r}\cr
&\qquad -\left(\frac{t}{r}\right)DB(y;U)t^2\overline{\nabla}f
-\frac{1}{t}\left(\frac{t}{r}\right)^2B(y)t^2(\overline{\opHess}(f)\circ\pi)(U)\cr
&\qquad+\frac{1}{t}\left(\frac{t}{r}\right)^2\left\langle U,\frac{y}{r}\right\rangle B(y)t^2\overline{\nabla} f
-\frac{1}{t}\left(\frac{t}{r}\right)^2\langle U, t^2\overline{\nabla}f\rangle B(y)\frac{y}{r}
+\Gamma(U,\nabla\hat{f}),\cr}
$$
and, bearing in mind that $B(y)y=y$ and $\langle y,\overline{\nabla}f\rangle = 0$, we obtain
$$\eqalign{
\Delta\hat{f} &= \frac{m}{t}\left(\frac{t}{r}\right)
-\left(\frac{t}{r}\right)\opTr(DB(y;\cdot)t^2\overline{\nabla}f)
-\frac{1}{t}\left(\frac{t}{r}\right)^2\opTr(B(y)t^2(\overline{\opHess}(f)\circ\pi))\cr
&\qquad + \frac{1}{t}\opTr(\Gamma(\cdot,tx))
-\left(\frac{t}{r}\right)\opTr(\Gamma(\cdot,B(y)t^2\overline{\nabla}f))\cr
&=\frac{1}{t}\opTr(\Phi_1) + \opTr(\Phi_2) + \frac{1}{t}\opTr(\Phi_3 t^2\overline{\opHess}(f)\circ\pi),\cr}
$$
where
$$
\Phi_1,\Phi_2,\Phi_3 \in \Bbb{R}_*[F]\otimes\overline{\Cal{Q}}^0_{*,*}[X,\overline{\nabla}F]
\otimes\overline{\Cal{R}}^0_{*,*}[X,\overline{\nabla}F][[T]]_{(2,1,2,1,2)}.
$$
Likewise,
$$\eqalign{
g(\nabla_{\nabla\hat{f}}\nabla\hat{f},\nabla\hat{f})
&=\frac{1}{t}\left(\frac{t}{r}\right)^3\langle B(y)t^2\overline{\nabla} f, t^2\overline{\nabla} f\rangle+\langle A(y)\Gamma(\nabla\hat{f},\nabla\hat{f}),\nabla\hat{f}\rangle\cr
&\qquad-\left(\frac{t}{r}\right)\langle A(y)DB(y;\nabla\hat{f})t^2\overline{\nabla} f\nabla\hat{f}\rangle\cr
&\qquad-\frac{1}{t}\left(\frac{t}{r}\right)^4\langle t^2(\overline{\opHess}\circ\pi)B(y)t^2\overline{\nabla}f,B(y)t^2\overline{\nabla}f\rangle.\cr}
\eqnum{\nexteqnno[SecondPartOfMeanCurvature]}
$$
However, for any symmetric bilinear form, $M_{ij}$, and for any vector $V^i$,
$$\eqalign{
\delta_{ij}\delta^{ip}M_{pq}B^q_r V^r B^j_s V^s
&= (\delta^{ip}\delta^{jq})M_{pq}(\delta_{jb}B^b_cV^c)(\delta_{ir}B^r_sV^s)\cr
&= (\delta^{ip}\delta^{jq})M_{pq}(B^b_j\delta_{bc}V^c)(B^r_i\delta_{rs}V^s),\cr}
$$
so that
$$
\langle M B(y)t^2\overline{\nabla}f,B(y)t^2\overline{\nabla}f\rangle = \langle M,\Psi\rangle,
$$
for some $\Psi\in\Bbb{R}_*[F]\otimes\overline{\Cal{Q}}^2_{*,*}[X,\overline{\nabla}F][[T]]_{(2,1,2)}$.
\par
Now, by \eqnref{GradientOfF}, $\nabla\hat{f}$ contains a term in $x$ that does not carry a factor of $t$. We need to show that this term in $x$ is not repeated in any non-trivial component of \eqnref{SecondPartOfMeanCurvature}. However, since $D_x x=\nabla_x x=0$, we have $\Gamma(x,x)=0$, so that, for all $U$,
$$
\langle A(y)\Gamma(x,x),U\rangle = 0.
$$
Next, since $g(y)(x,x)=1$ for all $y$, we obtain, for all vectors $U$,
$$
0 = g(y)(\nabla_U x,x)= g(y)(D_U x + \Gamma(U,x),x),
$$
so that
$$
\langle A(y)\Gamma(U,x),x\rangle = \langle A(y)\Gamma(x,U),x\rangle = -\langle A(y)D_U x,x\rangle = -\frac{1}{r}\langle A(y)\pi(U),x\rangle = 0.
$$
Finally, since $D_x\overline{\nabla}f=0$, and since $\langle B(y) t^2\overline{\nabla} f,x\rangle =0$ for all $y$, we have
$$
\langle DB(y,x)t^2\overline{\nabla}f,x\rangle = D_x\langle B(y)t^2\overline{\nabla}f,x\rangle = 0,
$$
and we conclude that the term in $x$ is not repeated in any non-trivial component of \eqnref{SecondPartOfMeanCurvature}, as desired. It follows that
$$
g(\nabla_{\nabla\hat{f}}\nabla\hat{f},\nabla\hat{f}) = \frac{1}{t}\Phi_4 + \Phi_5 + \langle\Phi_6,t^2(\overline{\opHess}\circ\pi)\rangle,
$$
where
$$\eqalign{
\Phi_4,\Phi_5 &\in \Bbb{R}_*[F]\otimes\overline{\Cal{Q}}^0_{*,*}[X,\overline{\nabla} F][[T]]_{(2,1,2)},\ \&\cr
\Phi_6 &\in \Bbb{R}_*[F]\otimes\overline{\Cal{Q}}^2_{*,*}[X,\overline{\nabla} F][[T]]_{(2,1,2)}.\cr}
$$
and the result for $Y=0$ now follows by Proposition \procref{TaylorSeriesOfNormN}. The general case follows by substituting for every generator $R_{i_1i_2i_3}{}^j{}_{;i_4...i_{k+3}}$ of $\Cal{R}$ its own Taylor series in $Y$ about $0$, as in the proof of Proposition \procref{TaylorSeriesOfNVariableCentre}. This completes the proof.\qed
\newhead{Asymptotic Expansions and Formal Solutions}[AsymptoticExpansionsAndFormalSolutions]
\newsubhead{Asymptotic Expansions}[AsymptoticExpansions]
In order to save on notation, which would otherwise quickly get out of hand, we shall no longer be so explicit about the definition of curvature polynomials, leaving the reader to infer how they are constructed in each case. We now reformulate the results of the previous sections in a manner that will allow us to construct formal solutions later on. To this end, we introduce the terminology of asymptotic expansions for functions defined near $t=0$. Thus, let $E$ be a finite-dimensional vector bundle over some finite-dimensional base $B$. Let $\phi:]0,\infty[\times E\rightarrow\Bbb{R}$ be a smooth function. Let $(\phi_k)$ be a sequence of smooth functions, where, for all $k$, $\phi_k:E^{\otimes k}\rightarrow\Bbb{R}$. For a formal power series $\xi_x(t)\sim\sum_{k=0}^\infty t^k\xi_{k,x}$ in $E$, we write
$$
\phi(t,\xi_x) \sim \sum_{k=0}^\infty t^k\phi_k(\xi_{0,x},...,\xi_{k,x})\eqnum{\nexteqnno[DefinitionOfAsymptoticExpansionPartI]}
$$
to mean that for all $N\geq 0$, there exists a smooth function $R_N:[0,\infty[\times E^{\otimes (N+1)}\rightarrow\Bbb{R}$ such that
$$
\phi\left(t,\sum_{k=0}^N t^k\xi_{x,k}\right) = \sum_{k=0}^Nt^k\phi_k(\xi_{0,x},...,\xi_{k,x}) + t^{N+1}R_N(t,\xi_{x,0},...,\xi_{x,N}).
\eqnum{\nexteqnno[DefinitionOfAsymptoticExpansionPartII]}
$$
It is of fundamental importance in our definition that the remainder term, $R_N$, be smooth also at $t=0$, as it would otherwise be of little use to us.
\proclaim{Proposition \nextprocno}
\noindent There exists a sequence $(P_k)$ of curvature polynomials such that for all formal power series $Y\sim\sum_{k=0}^\infty t^kY_k$ of vectors in $\Bbb{R}^{m+1}$ and $f_x\sim\sum_{k=0}^\infty t^k f_{k,x}$ of germs in $J^0S^m$, and for all vectors $U$,
$$
t\langle p(t,y,f_x),U\rangle \sim t\langle x,U\rangle + \sum_{k=0}^\infty t^k\langle P_k(f_{0,x},...,f_{k-3,x},Y_0,...,Y_{k-3}),U\rangle,
\eqnum{\nexteqnno[AsymptoticSeriesForP]}
$$
and $P_k=0$ for $k\leq 2$.
\endproclaim
\proclabel{AsymptoticSeriesForP}
\proof By Proposition \procref{TaylorSeriesOfLittleP}, there exists a smooth function $\tilde{P}$ such that $\langle p(t,Y,f_x),U\rangle = \langle x,U\rangle + t^2\langle\tilde{P}(t,Y,f_x),U\rangle$. Furthermore, since the coefficients of the Taylor series of $\tilde{P}$ in $t$ are all curvature polynomials, there exists a sequence $(\tilde{P}_k)$ of curvature polynomials such that
$$
\tilde{P}(t,Y,f_x) \sim \sum_{k=0}^\infty t^k\tilde{P}_k(f_{0,x},...,f_{k,x},Y_0,...,Y_k).
$$
It follows that
$$
t\langle p(t,Y,f_x),U\rangle \sim t\langle x,U\rangle + \sum_{k=3}^\infty t^k\langle\tilde{P}_{k-3}(f_{0,x},...,f_{k-3,x},Y_0,...,Y_{k-3}),U\rangle,
$$
as desired.\qed
\proclaim{Proposition \nextprocno}
\noindent There exists a sequence $(Q_k)$ of curvature polynomials such that for all formal power series $Y\sim\sum_{k=0}^\infty t^kY_k$ and $V\sim\sum_{k=0}^\infty t^k V_k$ of vectors in $\Bbb{R}^{m+1}$ and $f_x\sim\sum_{k=0}^\infty t^k f_{k,x}$ of germs in $J^0S^m$,
$$
t\langle q(t,Y,f_x),V\rangle \sim \sum_{k=0}^\infty t^k\left(\langle x,V_{k-2}\rangle + Q_k(f_{0,x},...,f_{k-4,x},Y_0,...,Y_{k-4},V_0,...,V_{k-4})\right),
\eqnum{\nexteqnno[AsymptoticSeriesForQ]}
$$
where $Q_k=0$ for $k\leq 3$.
\endproclaim
\proclabel{AsymptoticSeriesForQ}
\proof By Proposition \procref{TaylorSeriesOfLittleQ}, there exists a smooth function $\tilde{Q}$ such that $\langle q(t,Y,f_x),V\rangle=t\langle x,V\rangle + t^3\langle\tilde{Q}(t,Y,f_x),V\rangle$. Furthermore, since the coefficients of the Taylor series of $\tilde{Q}$ in $t$ are all curvature polynomials, there exists a sequence $(\tilde{Q}_k)$ of curvature polynomials such that
$$
\langle\tilde{Q}(t,Y,f_x),V\rangle \sim \sum_{k=0}^\infty t^k \tilde{Q}_k(f_{0,x},...,f_{k,x},Y_{0},...,Y_{k},V_0,...,V_k).
$$
It follows that,
$$
t\langle q(t,Y,f_x),V\rangle \sim \sum_{k=2}^\infty t^k\langle x,V_{k-2}\rangle +
\sum_{k=4}^\infty \tilde{Q}_{k-4}(f_{0,x},...,f_{k-4,x},Y_0,...,Y_{k-4},V_0,...,V_{k-4}),
$$
as desired.\qed
\proclaim{Proposition \nextprocno}
\noindent There exists a sequence $(R_k)$ of curvature polynomials such that for all formal power series $Y\sim\sum_{k=0}^\infty t^kY_k$ of vectors in $\Bbb{R}^{m+1}$ and $f_x\sim\sum_{k=0}^\infty t^kf_{k,x}$ and $g_x\sim\sum_{k=0}^\infty t^kg_{k,x}$ of germs in $J^0S^m$,
$$
tr(t,Y,f_x)g_x \sim \sum_{k=0}^\infty t^k\left(t^4 g_{k,x} + R_k(Y_0,...,Y_{k-4},f_{0,x},...,f_{k-4,x},t^4g_{0,x},...,t^4g_{k-4,x})\right),
\eqnum{\nexteqnno[AsymptoticSeriesForR]}
$$
where $R_k=0$ for $k\leq 4$.
\endproclaim
\proclabel{AsymptoticSeriesForR}
\proof By Proposition \procref{TaylorSeriesOfLittleR}, there exists a smooth function $\tilde{R}$ such that $r(t,Y,f_x)=t^3+t^7\tilde{R}(t,Y,f_x)$. Furthermore, since the coefficients of the Taylor series of $\tilde{R}$ in $t$ are all curvature polynomials, there exists a sequence $(\tilde{R}_k)$ of curvature polynomials such that
$$
\tilde{R}(t,Y,f_x) \sim \sum_{k=0}^\infty t^k\tilde{R}_k(Y_0,...,Y_k,f_{0,x},...,f_{k,x}).
$$
Thus
$$
t^4\tilde{R}(t,Y,f_x)g_x \sim \sum_{k=0}^\infty t^k\sum_{l=0}^k\tilde{R}_l(Y_0,...,Y_l,f_{0,x},...,f_{l,x})(t^4g_{k-l,x}).
$$
It follows that
$$
tR(t,Y,f_x)g_x \sim \sum_{k=0}^\infty t^k(t^4g_{k,x}) + \sum_{k=4}^\infty t^k\sum_{l=0}^{k-4}\tilde{R}_l(Y_0,...,Y_l,f_{0,x},...,f_{l,x})(t^4g_{k-l-4,x}),
$$
as desired.\qed
\proclaim{Proposition \nextprocno}
\noindent There exists a sequence $(H_k)$ of curvature polynomials such that for all formal power series $Y\sim\sum_{k=0}^\infty t^kY_k$ of vectors in $\Bbb{R}^{m+1}$ and $f_x\sim\sum_{k=0}^\infty t^kf_{k,x}$ of germs in $J^2S^m$,
$$\eqalign{
\frac{1}{t}\left(H(t,Y,f)-\frac{1}{t}\right) &\sim\sum_{k=0}^\infty t^k\left( -\frac{1}{n}(n+\Delta)f_{k,x} - \frac{1}{4}\opRic_{pq;rs}x^px^qx^rY^s_{k-2}\right.\cr
&\qquad \left.-\frac{1}{3}\opRic_{pq;r}x^px^qY^r_{k-1} + H_k(Y_0,...,Y_{k-3},f_{0,x},...,f_{k-2,x})\right),\cr}
\eqnum{\nexteqnno[AsymptoticExpansionOfH]}
$$
where, by convention, $Y_k=0$ for $k<0$. Furthermore,
$$\eqalign{
H_0 &= -\frac{1}{3}\opRic_{pq}x^px^q,\ \&\hfill\cr
H_1 &= -\frac{1}{4}\opRic_{pq;r}x^px^qx^r.\hfill\cr}
\eqnum{\nexteqnno[FirstTermsInAsymptoticExpansionOfH]}
$$
\endproclaim
\proclabel{AsymptoticSeriesForH}
\proof Consider the formula \eqnref{FirstAsymptoticFormulaForDisplacedCurvature} for $H$. Trivially,
$$\eqalign{
\frac{1}{n}(n+\overline{\Delta})f_x&\sim \sum_{k=0}^\infty t^k\frac{1}{n}(n+\overline{\Delta})f_{k,x},\cr
\frac{t}{3}\opRic_{pq;r}x^px^qY^r&\sim \sum_{k=1}^\infty \frac{t^k}{3}\opRic_{pq;r}x^px^qY^r_{k-1},\ \&\cr
\frac{t^2}{4}\opRic_{pq;rs}x^px^qx^rY^s&\sim \sum_{k=2}^\infty \frac{t^k}{4}\opRic_{pq;r}x^px^qx^rY^r_{k-2}.\cr}
$$
Since $F$ is a curvature polynomial, there exists a sequence $(F_k)$ of curvature polynomials such that
$$
F(f_x) \sim \sum_{k=0}^\infty t^k F_k(f_{0,x},...,f_{k,x}).
$$
In particular
$$
t^2F(f_x) \sim \sum_{k=2}^\infty t^k F_{k-2}(f_{0,x},...,f_{k-2,x}).
$$
Finally, denote the remainder term in \eqnref{FirstAsymptoticFormulaForDisplacedCurvature} by $t^5G(t,Y,f)$. Since every coefficient in the Taylor series of $G$ in $t$ about $0$ is a curvature polynomial, there exists a sequence $(G_k)$ of curvature polynomials such that
$$
G(t,Y,f) \sim \sum_{k=0}^\infty t^k G_k(Y_0,...,Y_k,f_{0,x},...,f_{k,x}).
$$
In particular
$$
t^3G(t,Y,f) \sim \sum_{k=3}^\infty t^k G_{k-3}(Y_0,...,Y_{k-3},f_{0,x},...,f_{k-3,x}),
$$
and the result follows upon combining these terms.\qed
\newsubhead{Flows of Surfaces}[FlowsOfSurfaces]
We now extend our framework to the time-dependent case which interests us. Thus, let $M$ be an $(m+1)$-dimensional Riemannian manifold with metric $g$, let $R$ be its Riemann curvature tensor, let $S$ be its scalar curvature function, and suppose that $S$ is of Morse-Smale type. Let $\gamma:\Bbb{R}\rightarrow M$ be a complete integral curve of $-\nabla S$ with relatively compact image. In particular (c.f. \cite{Schwarz}), $\gamma(t)$ converges exponentially to critical points of $S$ as $t$ tends to $\pm\infty$, and its derivatives to all orders decay exponentially at infinity.
\par
For convenience, we suppose that $M$ has unit injectivity radius. We identify the bundle $\gamma^*TM$ with the product bundle $\Bbb{R}\times\Bbb{R}^{m+1}$ via parallel transport. For all $t\in\Bbb{R}$, define the metric $g_t$ over $\Bbb{R}^{m+1}$ by $g_t:=\opExp_{\gamma(t)}^*g$, where $\opExp_{\gamma(t)}$ here denotes the exponential map of $M$ about the point $\gamma(t)$. In particular, for all $t$, the metric $g_t$ is of the type introduced in Section \subheadref{CurvatureTensors}. Furthermore, the family $(g_t)$ converges exponentially in the $C^\infty_\oploc$ sense to metrics $g_{\pm\infty}$ as $t$ tends to $\pm\infty$ and its time derivatives to all orders also decay exponentially at infinity.
\par
As in Section \subheadref{TheExponentialMapAndParallelTransport}, for all $t\in\Bbb{R}$, let $\opExp_t:\Omega_2\rightarrow\Bbb{R}^{m+1}$ be the exponential map of $g_t$. That is, for all $(x,y)\in\Omega_2$, the curve $s\mapsto \opExp_t(x,sy)$ is the unique geodesic with respect to $g_t$ leaving the point $x$ in the direction of the vector $y$. For $s>0$, and for bounded functions $Y\in C^0(\Bbb{R},\Bbb{R}^{m+1})$ and $f\in C^0(\Bbb{R}\times S^m)$, define $e(s,Y,f):\Bbb{R}\times S^m\rightarrow\Bbb{R}^{m+1}$ by
$$
e(s,Y,f)(t,x):=\opExp_t(sY(t),s(1+s^2f(t,x))x).\eqnum{\nexteqnno[FamilyOfImmersions]}
$$
Heuristically, $e(s,Y,f)$ is a continuous family of immersed spheres all of radius approximately $s$, with centres displaced by the function $Y$. Composing with $\opExp_{\gamma(s)}$ then yields a continuous family of small immersed spheres in $M$ which move along the geodesic $\gamma$. We will show that for sufficiently small $s$ and for correct choices of $Y$ and $f$, this family yields a forced mean curvature flow of immersed spheres in $M$ with forcing term $1/s$.
\par
For all $k$, let $J^k(\Bbb{R},\Bbb{R}^{m+1})$ denote the bundle of $k$-jets over $\Bbb{R}$ taking values in $\Bbb{R}^{m+1}$. For all $(k,l)$, let $J^{k,l}(\Bbb{R}\times S^m,\Bbb{R})$ denote the bundle of $(k,l)$-jets over $\Bbb{R}\times S^m$ taking values in $\Bbb{R}$, that is, the bundle of $\Bbb{R}$-valued jets that are of order at most $k$ in $\Bbb{R}$ and at most $l$ in $S^m$. Observe that $J^{k,l}(\Bbb{R}\times S^m,\Bbb{R})$ is actually also a bundle over $\Bbb{R}$ and we denote by $J:=J^{k,l}$ its fibrewise cartesian product with $J^k(\Bbb{R},\Bbb{R}^{m+1})$. In other words, an element of $J^{k,l}$ is a pair $(Y_t,f_{t,x})$ where $Y_t$ is the jet of an $\Bbb{R}^{m+1}$-valued function over $\Bbb{R}$ at the point $t$, and $f_{t,x}$ is the jet of an $\Bbb{R}$-valued function over $\Bbb{R}\times S^m$ at the point $(t,x)$.
\par
Define the functions $N:]0,\infty[\times J\rightarrow S^m$ and $H:]0,\infty[\times J\rightarrow\Bbb{R}$ such that for all $s\in]0,\infty[$ and for all $(Y_t,f_{t,x})\in J$, $N(s,Y_t,f_{t,x})$ and $H(s,Y_t,f_{t,x})$ are respectively the outward-pointing unit normal of the immersion $e(s,Y,f)(t,\cdot)$ at the point $e(s,Y,f)(t,x)$ and its mean curvature at that point, both with respect to the metric $g_t$. Define $V:]0,\infty[\times J\rightarrow\Bbb{R}^{m+1}$ by
$$
V(s,(Y_t,f_{t,x})) := \partial_r \opExp_{\gamma(t)}^{-1}(\opExp_{\gamma(r)}(e(s,Y,f)(t+r,x)))|_{r=0}.\eqnum{\nexteqnno[VariationOfImmersions]}
$$
Heuristically, this vector field measures the variation of the immersion $e(s,Y,f)$ at the point $e(s,Y,f)(t,x)$ as we move along the flow. Finally, define $\Phi:]0,\infty[\times J\rightarrow\Bbb{R}$ by
$$
\Phi(s,(Y_t,f_{t,x})) := \frac{1}{s}\left(H(s,(Y_t,f_{t,x})) - \frac{1}{s}\right) + s\langle V(s,(Y_t,f_{t,x})),N(s,(Y_t,f_{t,x}))\rangle.\eqnum{\nexteqnno[ForcedMCFEquation]}
$$
For all $s$, $\Phi(s,\cdot)$ is the {\bf forced mean curvature flow operator} (with forcing term $1/s$). In particular, it is a quasi-linear parabolic partial differential operator whose zeroes are (reparametrised) forced mean curvature flows with forcing term $1/s$.
\proclaim{Proposition \nextprocno}
\noindent There exists a sequence $(\Phi_k)$ of curvature polynomials such that for all formal power series $(Y_t,f_{t,x})\sim \sum_{k=0}^\infty s^k(Y_t,f_{t,x})$ of germs in $J$,
$$\multiline{
\Phi(s,Y_t,f_{t,x})\sim \sum_{k=0}^\infty s^k\left[\left\langle \left(\frac{\partial}{\partial t}
+\frac{(m+1)}{2(m+3)}\opHess(S)\right)Y_{k-2,t},x\right\rangle\right.\cr
\qquad\qquad+\left(s^4\frac{\partial}{\partial t} + \frac{1}{m}(m+\overline{\Delta})\right)f_{k,x,t}\cr \qquad\qquad+\left(\frac{1}{4}\opRic_{t,ab;cd}x^ax^bx^cY^d_{k-2,t} - \frac{(m+1)}{2(m+3)}S_{t,;ab}x^aY^b_{k-2,t}\right)\cr
\qquad\qquad-\frac{1}{3}\opRic_{t,ab;c}x^ax^bY^c_{k-1,t}\cr
\qquad\qquad+\Phi_k(f_{0,x,t},...,f_{k-2,x,t},s^4\dot{f}_{0,x,t},...,s^4\dot{f}_{k-4,x,t},Y_{0,t},...,Y_{k-3,t},\dot{Y}_{0,t},...,\dot{Y}_{k-4,t})\bigg],\cr}
\eqnum{\nexteqnno[AsymptoticExpansionOfPhi]}
$$
where $\opRic_t$ and $S_t$ denote respectively the Ricci and Scalar curvatures of $M$ at the point $\gamma(t)$, and, by convention, $Y_k=0$ for $k<0$. Furthermore, the curvature polynomials $\Phi_0$ and $\Phi_1$ are given by
$$\eqalign{
\Phi_0 &= -\frac{1}{3}\opRic_{ab}x^ax^b,\cr
\Phi_1 &= -\frac{1}{4}\opRic_{ab;c}x^ax^bx^c+\frac{(m+1)}{2(m+3)}S_{;a}x^a.\cr}
\eqnum{\nexteqnno[FirstTermsInAsymptoticSeriesOfPhi]}
$$
\endproclaim
\proclabel{AsymptoticExpansionOfPhi}
\remark Importantly, since they are curvature polynomials, the functions $(\Phi_k)$ vary with $t$ only insofar as the curvature tensor itself, along with its derivatives, vary with $t$, and the same can also be said for the remainder terms in the asymptotic series. In particular, since the flow $\gamma$ has relatively compact image in $M$, the derivatives of all these functions to all orders are uniformly bounded independent of $s$ and $t$.
\medskip
\remark Observe that, as in Proposition \procref{AsymptoticSeriesForR}, in every remainder term of this asymptotic series, the term $\dot{f}_k$ only ever appears accompanied by the factor $s^4$.
\medskip
\proof Indeed
$$
V(t,(Y_t,f_{t,x})) = P(t,(Y_t,f_{t,x}))\dot{\gamma} + Q(t,(Y_t,f_{t,x}))\dot{Y}_t + R(t,(Y_t,f_{t,x}))\dot{f}_{t,x}.
$$
Furthermore, since $\gamma$ is a gradient flow of $S$,
$$
\dot{\gamma} = -\frac{(m+1)}{2(m+3)}S_{;a}x^a.
$$
The result now follows by Propositions \procref{AsymptoticSeriesForP}, \procref{AsymptoticSeriesForQ}, \procref{AsymptoticSeriesForR} and \procref{AsymptoticSeriesForH}.\qed
\newsubhead{Parabolic Operators I - The Finite Dimensional Case}[ParabolicOperatorsI]
We first aim to determine formal solutions of the equation $\Phi(s,Y,f)=0$ for small values of $s$. To this end, we introduce the following functional analytic framework. For a finite-dimensional vector space, $E$, and for $\alpha\in]0,1]$, define the {\bf H\"older seminorm} of order $\alpha$ over $C^0(\Bbb{R},E)$ by
$$
[f]_\alpha:=\msup_{0<\left|s-t\right|\leq 1}\frac{\|f(s)-f(t)\|}{\left|s-t\right|^\alpha}.\eqnum{\nexteqnno[DefinitionOfHolderSeminorm]}
$$
For all $k$ and for all $\alpha\in]0,1]$, define the {\bf H\"older norm} of order $(k,\alpha)$ over $C^k(\Bbb{R},E)$ by
$$
\|f\|_{k,\alpha} := \sum_{i=0}^k\|\partial_t^i f\|_0 + [\partial_t^k f]_\alpha,\eqnum{\nexteqnno[DefinitionOfHolderNorm]}
$$
where $\|\cdot\|_0$ denotes the uniform norm. For all $(k,\alpha)$, define the {\bf H\"older space} of order $(k,\alpha)$ by
$$
C^{k,\alpha}(\Bbb{R},E) := \left\{ f\in C^k(\Bbb{R},E)\ |\ \|f\|_{k,\alpha} < \infty \right\}.\eqnum{\nexteqnno[DefinitionOfHolderSpace]}
$$
Recall that $C^{k,\alpha}$ furnished with the norm $\|\cdot\|_{k,\alpha}$ constitutes a Banach space.
\par
Define the operator $P:C^{1,\alpha}(\Bbb{R},\Bbb{R}^{m+1})\rightarrow C^{0,\alpha}(\Bbb{R},\Bbb{R}^{m+1)}$ by
$$
PY = \left(\frac{\partial}{\partial t}
+\frac{(m+1)}{2(m+3)}\opHess(S)\right)Y.\eqnum{\nexteqnno[DefinitionOfOperatorP]}
$$
Observe that this operator corresponds to the first summand in the asymptotic expansion \eqnref{AsymptoticExpansionOfPhi} of $\Phi$. Furthermore, since $S$ is of Morse-Smale type, $P$ is Fredholm and surjective. In addition, since every function in $\opKer(P)$ decays exponentially at infinity (c.f. \cite{Schwarz}), the $L^2$ orthogonal complement, $\opKer(P)^\perp$, of $\opKer(P)$ in $C^{1,\alpha}(\Bbb{R},\Bbb{R}^{m+1})$ is well-defined. The restriction of $P$ to $\opKer(P)^\perp$ is invertible, and we denote its inverse by $G$.
\par
We will also be interested in families of constant coefficient parabolic operators over $C^1(\Bbb{R},E)$. Thus, for an invertible linear map $A:E\rightarrow E$, which, for convenience, we take to be symmetric with respect to some fixed metric over $E$, and for $\epsilon>0$, define $P_\epsilon:C^1(\Bbb{R},E)\rightarrow C^0(\Bbb{R},E)$ by
$$
P_\epsilon f := (\epsilon\partial_t - A)f.\eqnum{\nexteqnno[DefinitionOfPEpsilon]}
$$
It follows from the invertibility of $A$ that $P_\epsilon$ as also invertible. In fact, its Green's operator, which we denote by $G_\epsilon$, is given by
$$
G_\epsilon f(t) = -\frac{1}{\epsilon}\int_{-\infty}^t e^{-\frac{1}{\epsilon}(t-s)A^+}f(s)ds + \frac{1}{\epsilon}\int_t^\infty e^{-\frac{1}{\epsilon}(t-s)A^-}f(s)ds.\eqnum{\nexteqnno[DefinitionOfGEpsilon]}
$$
where $A^+$ (resp. $A^-$) denotes the composition of $A$ with the orthogonal projection onto the direct sum of its eigenspaces of positive (resp. negative) eigenvalue. In order to obtain uniform estimates on the operator norm of $G_\epsilon$, it is useful to introduce a weighting factor into the H\"older norm. Thus, for all $(k,\alpha)$ and for all $\epsilon>0$, define the {\bf weighted H\"older norm} of order $(k,\alpha)$ and weight $\epsilon$ by
$$
\|f\|_{k,\alpha,\epsilon} := \sum_{i=0}^k\epsilon^i\|\partial_t^if\|_0 + \epsilon^k[\partial_t^if]_\alpha\eqnum{\nexteqnno[FirstWeightedHolderNorm]}
$$
Observe that, for all $\epsilon$, the norm $\|\cdot\|_{k,\alpha,\epsilon}$ is uniformly equivalent to the norm $\|\cdot\|_{k,\alpha}$, so that $C^{k,\alpha}(\Bbb{R}\times E)$ is also a Banach space with respect to every weighted H\"older norm.
\proclaim{Proposition \nextprocno}
\noindent There exists $B>0$, which only depends on the matrix $A$, such that for all $\epsilon>0$, and for all $f\in C^{0,\alpha}(\Bbb{R},E)$,
$$
\|G_\epsilon f\|_{1,\alpha,\epsilon} \leq B\|f\|_{0,\alpha}
\eqnum{\nexteqnno[OperatorNormOfGreensOperator]}
$$
\endproclaim
\proclabel{NormOfInverseOfFiniteDimensionalParabolicOperator}
\proof Since both $P_\epsilon$ and $G_\epsilon$ preserve the eigenspaces of $A$, we may suppose that $E=\Bbb{R}$ and that $A=\lambda>0$. Thus
$$
G_\epsilon f(t) = -\frac{1}{\epsilon}\int_{-\infty}^t e^{-\frac{\lambda(t-s)}{\epsilon}}f(s)ds
=-\frac{1}{\epsilon}\int_0^\infty e^{-\frac{\lambda s}{\epsilon}}f(t-s)ds.
$$
Now fix $f\in C^{0,\alpha}(\Bbb{R},\Bbb{R})$. For all $t$,
$$
\left|G_\epsilon f(t)\right| \leq \frac{1}{\epsilon}\int_0^\infty e^{-\frac{\lambda}{\epsilon s}}\left|f(s)\right|ds \leq \frac{1}{\lambda}\|f\|_0,
$$
and taking the supremum over all $t$ yields $\|G_\epsilon f\|_0\leq \frac{1}{\lambda}\|f\|_0$. Likewise, for all $0<\left|t-t'\right|\leq 1$,
$$
\left|G_\epsilon f(t) - G_\epsilon f(t')\right| \leq \frac{1}{\epsilon}\int_0^\infty e^{-\frac{\lambda}{\epsilon}s}\left|f(t-s) - f(t'-s)\right| \leq \frac{1}{\lambda}\left|t-t'\right|^\alpha [f]_\alpha.
$$
Dividing both sides by $\left|t-t'\right|^\alpha$, and taking the supremum over all $t$ yields $[G_\epsilon f]_\alpha\leq\frac{1}{\lambda}[f]_\alpha$. Combining these relations yields $\|G_\epsilon f\|_{0,\alpha}\leq \|A\|^{-1}\|f\|_{0,\alpha}$. Finally, by definition of $G_\epsilon$, $\epsilon\partial_t G_\epsilon f = \lambda G_\epsilon f + f$, so that $\epsilon\|\partial_t G_\epsilon f\|_{0,\alpha}\leq \lambda\|G_\epsilon f\|_{0,\alpha} + \|f\|_{0,\alpha} \leq 2\|f\|_{0,\alpha}$. This completes the proof.\qed
\newsubhead{Parabolic Operators II - the Infinite-Dimensional Case}[ParabolicOperatorsII]
For all $\alpha\in]0,1]$, define the {\bf H\"older seminorms} of order $\alpha$ over $C^0(\Bbb{R}\times\Bbb{S}^m)$ by
$$\eqalign{
[f]_{x,\alpha} \hfill&:= \msup_{t,x\neq y} \frac{\left|f(t,x) - f(t,y)\right|}{\|x-y\|^\alpha},\cr
[f]_{t,\alpha} \hfill&:= \msup_{x,0\left|t-s\right|\leq 1} \frac{\left|f(s,x) - f(t,x)\right|}{\left|s-t\right|^\alpha}.\cr}\eqnum{\nexteqnno[DefinitionOfInhomogeneousHolderSeminorms]}
$$
For all $k\in\Bbb{N}$, let $C^k_\opin(\Bbb{R}\times S^m)$ be the set of all functions $f:\Bbb{R}\times S^m\rightarrow\Bbb{R}$ which are continuously differentiable $i$ times in the $x$ direction and $j$ times in the $t$ direction for all $i+2j\leq 2k$. For all $k\in\Bbb{N}$ and for all $\alpha\in ]0,1/2]$, define the {\bf inhomogeneous H\"older norm} of order $(k,\alpha)$ over $C^k_\opin(\Bbb{R}\times S^m)$ by
$$
\|f\|_{k,\alpha,\opin} := \sum_{i+2j\leq 2k}\|D_x^i D_t^j f\|_0 + \sum_{i+2j=2k}[D_x^iD_t^j f]_{x,2\alpha} + \sum_{i+2j=2k}[D_x^iD_t^j f]_{t,\alpha}.\eqnum{\nexteqnno[DefinitionOfInhomogeneousHolderNorm]}
$$
For all $k,\alpha$, define the {\bf inhomogeneous H\"older space} of order $(k,\alpha)$ by
$$
C^{k,\alpha}_\opin(\Bbb{R}\times S^m) := \left\{ f\in C^k_\opin(\Bbb{R}\times S^m)\ |\ \|f\|_{k,\alpha,\opin}<\infty\right\}.\eqnum{\nexteqnno[DefinitionOfInhomogeneousHolderSpace]}
$$
Recall that $C^{k,\alpha}_\opin(\Bbb{R}\times S^m)$ furnished with the norm $\|\cdot\|_{k,\alpha,\opin}$ constitutes a Banach space. More generally, for all $(k,\alpha)$ and for all $\epsilon>0$, define the {\bf weighted inhomogeneous H\"older norm} of order $(k,\alpha)$ and weight $\epsilon$ over $C^k_\opin(\Bbb{R}\times S^m)$ by
$$
\|f\|_{k,\alpha,\opin,\epsilon} := \sum_{i+2j\leq 2k}\epsilon^j\|D_x^i D_t^j f\|_0 + \sum_{i+2j=2k}\epsilon^j[D_x^iD_t^j f]_{x,2\alpha} + \sum_{i+2j=2k}\epsilon^j[D_x^iD_t^j f]_{t,\alpha}.\eqnum{\nexteqnno[DefinitionOfWeightedHolderNorm]}
$$
For all $\epsilon>0$, the norm $\|\cdot\|_{k,\alpha,\opin,\epsilon}$ is uniformly equivalent to the norm $\|\cdot\|_{k,\alpha,\opin}$ and it follows that $C^{k,\alpha}_\opin(\Bbb{R}\times S^m)$ is also a Banach space with respect to every weighted inhomogeneous H\"older norm.
\par
For all $s>0$, define the operator $Q_s:C^{1,\alpha}_\opin(\Bbb{R}\times S^m)\rightarrow C^{0,\alpha}_\opin(\Bbb{R}\times S^m)$ by
$$
Q_s f:=\left(s^4\frac{\partial}{\partial t} + \frac{1}{m}\left(m+\overline{\Delta}\right)\right)f,
\eqnum{\nexteqnno[DefinitionOfOperatorQ]}
$$
where, as in Section \headref{TaylorSeriesOfFunctionsDerivedFromImmersions}, $\overline{\Delta}$ denotes the Laplacian of the standard metric over $S^m$. Observe that this operator corresponds to the second summand in the asymptotic expansion \eqnref{AsymptoticExpansionOfPhi} of $\Phi$. Furthermore, the operator $(m+\overline{\Delta})$ defines a self-adjoint operator over $L^2(S^m)$ with kernel $\Cal{H}_1$, the space of restrictions to $S^m$ of linear functions over $\Bbb{R}^{m+1}$. In particular, $(m+\overline{\Delta})$ restricts to an invertible mapping of $\Cal{H}_1^\perp$ to itself. With this in mind, for all $k$ and for all $\alpha$, we define
$$
\hat{C}^{k,\alpha}_\opin(\Bbb{R}\times S^m)
:=\left\{ f\in C^{k,\alpha}_\opin(\Bbb{R}\times S^m)\ |\ \int_{S^m}f(t,x)x^i\opdVol = 0\ \forall 1\leq i\leq m+1\right\},
\eqnum{\nexteqnno[FunctionsInOrthogonalComplement]}
$$
and it follows from the classical theory of parabolic operators that, for all $s$, $Q_s$ restricts to an invertible mapping from $\hat{C}^{1,\alpha}_\opin(\Bbb{R}\times S^m)$ into $\hat{C}^{0,\alpha}_\opin(\Bbb{R}\times S^m)$. Uniform norm estimates for Green's operators in the infinite-dimensional setting differ significantly from those obtained in the finite-dimensional setting. Indeed,
\proclaim{Lemma \nextprocno}
\noindent There exists $B>0$ such that for all $s\leq 1$ and for all $f\in\hat{C}^{0,\alpha}_\opin(\Bbb{R}\times S^m)$
$$
\|H_sf\|_{1,\alpha,\opin,s^4} \leq Bs^{-4\alpha}\|f\|_{0,\alpha,\opin}.\eqnum{\nexteqnno[InfiniteDimensionalFactor]}
$$
\endproclaim
\proclabel{InfiniteDimensionalFactor}
\remark Although it may appear that this weaker estimate is merely a consequence of the naive approach to the proof, the study of solutions of the heat equation in euclidean space appears to indicate that it is probably optimal.
\medskip
\remark Alternatively, it may appear that this weaker estimate arises from the unusual definition \eqnref{DefinitionOfWeightedHolderNorm} of the weighted inhomogenous H\"older norm. Indeed, it would surely have made more sense to have multiplied the third summand of \eqnref{DefinitionOfWeightedHolderNorm} by a factor of $\epsilon^\alpha$, for in this case the factor of $s^{-4\alpha}$ would not have appeared in \eqnref{InfiniteDimensionalFactor}. However, we have chosen the above definition so that the operator $s^4\partial_t$ has unit norm with respect to the norms $\|\cdot\|_{1,\alpha,\opin,s^4}$ and $\|\cdot\|_{0,\alpha,\opin}$, as this ensures that other factors of $s^{-4\alpha}$ do not enter into our reasoning in places where they would be more of a technical nuisance.
\medskip
\proof For all $s>0$, define the isomorphism $D_s$ of $C^{k,\alpha}_\opin(\Bbb{R}\times S^m)$ by $D_sf(t,x)=f(s^4t,x)$. For all $s\leq 1$, and for all $f\in C^{0,\alpha}_\opin(\Bbb{R}\times S^m)$, $\|D_sf\|_{0,\alpha,\opin}\leq \|f\|_{0,\alpha,\opin}$. On the other hand, for all $s\leq 1$ and for all $f\in C^{1,\alpha}_\opin(\Bbb{R}\times S^m)$, $\|D_s^{-1}f\|_{1,\alpha,\opin,s^4}\leq t^{-4\alpha}\|f\|_{1,\alpha,\opin}$. However, for all $s$, $Q_s=D_s^{-1}Q_1D_s$. The result follows.\qed
\medskip
Observe that $\Cal{H}_1$ is really the space of eigenfunctions of $\overline{\Delta}$ of eigenvalue $m$. More generally, the decomposition of $L^2(S^m)$ into eigenspaces of $\overline{\Delta}$ actually yields better estimates for $\|H_sf\|_{1,\alpha,\opin,s^4}$ in the case where $f(t,\cdot)$ is the restriction to $S^m$ of an $s$-dependent polynomial function of bounded order. Indeed, for all $l$, let $\Cal{H}_l\subseteq L^2(S^m)$ be the space of {\bf spherical harmonics} of order $l$ over $S^m$, that is, the space of eigenfunctions of the operator $\overline{\Delta}$ with eigenvalue $l(m+l-1)$. Recall that, for all $l$, $\Cal{H}_l$ is the restriction to $S^m$ of the space of homogeneous harmonic polynomials of order $l$ over $\Bbb{R}^{m+1}$. In particular, any polynomial of order $l$ over $\Bbb{R}^{m+1}$ restricts to an element of $\Cal{H}_0\oplus...\oplus\Cal{H}_l$ over $S_m$. Now define
$$
\hat{\Cal{H}}_{l} := \oplus_{i=0,i\neq 1}^l\Cal{H}_i.\eqnum{\nexteqnno[NonLinearSphericalHarmonics]}
$$
Observe that $\hat{\Cal{H}}_{l}$ is contained in $\Cal{H}_1^\perp$ for all $l$. Furthermore, for all $l$ and for all $(k,\alpha)$, $C^{k,\alpha}(\Bbb{R},\hat{\Cal{H}}_{l})$ naturally identifies with a subspace of $\hat{C}^{k,\alpha}_\opin(\Bbb{R}\times S^m)$. In particular, for all $s$, $Q_s$ restricts to a mapping from $C^{1,\alpha}(\Bbb{R},\hat{\Cal{H}}_{l})$ to $C^{0,\alpha}(\Bbb{R},\hat{\Cal{H}}_{l})$. Furthermore, this restriction is invertible for all $s$, and Proposition \procref{NormOfInverseOfFiniteDimensionalParabolicOperator} now yields
\proclaim{Proposition \nextprocno}
\noindent For all $l\in\Bbb{N}$, there exists $B_l>0$ such that for all $f\in C^{0,\alpha}(\Bbb{R},\hat{\Cal{H}}_{\leq l})$ and for all $\epsilon$,
$$
\|H_s f\|_{1,\alpha,\opin,s^4} \leq B_l\|f\|_{0,\alpha,\opin}.
$$
\endproclaim
\proclabel{NormOfRestrictionOfGreensOperator}
\newsubhead{More on Spherical Harmonics}[MoreOnSphericalHarmonics]
A tensor $T^{i_1...i_k}$ is said to be {\bf isotropic} whenever
$$
A^{i_1}_{j_1}...A^{i_k}_{j_k}T^{j_1...j_k} = T^{i_1...i_k},
\eqnum{\nexteqnno[IsotropicTensors]}
$$
for all $i_1,...,i_k$ and for every special-orthogonal matrix $A$. Given two symmetric tensors $T_1^{i_1...i_k}$ and $T_2^{i_1...i_l}$, their symmetric product is given by
$$
(T_1\odot T_2)^{i_1...i_{k+l}} := \sum_{\sigma\in\tilde{\Sigma}_{k,l}}
T_1^{i_{\sigma(1)}...1_{\sigma(k)}}T_2^{i_{\sigma(k+1)}...1_{\sigma(k+l)}},\eqnum{\nexteqnno[SymmetricProduct]}
$$
where $\tilde{\Sigma}_{k,l}$ denotes the set of permutations of the set $\left\{1,...,k+l\right\}$ such that $\sigma(1)<...<\sigma(k)$ and $\sigma(k+1)<...<\sigma(k+l)$. Let $\delta$ be as in Section \subheadref{CurvatureTensors}. In particular, $\delta$ is symmetric and isotropic. Furthermore, for all $k$, its $k$'th symmetric power, $\delta^{\odot k}$, is also a symmetric and isotropic tensor. In fact, up to rescaling, these are the only ones.
\proclaim{Lemma \nextprocno}
\noindent The space of symmetric, isotropic tensors of order $k$ is $1$-dimensional when $k$ is even, and $0$-dimensional when $k$ is odd.
\endproclaim
\proclabel{SymmetricIsotropicTensors}
\proof Indeed, the space of symmetric tensors of order $k$ is isomorphic to the space of homogeneous polynomials of the same order. However, since an $\opSO(m+1)$-invariant polynomial is constant over every sphere centred on the origin, it is determined by its restriction to any straight line passing through the origin, and when, in addition, this polynomial is homogeneous, it is determined by its value at a single point. It follows that this space has dimension at most $1$. Now observe that the restriction of a homogeneous polynomial to a straight line through the origin is even when its order is even, and odd when its order is odd. However, by $\opSO(m+1)$-invariance again, the restrictions of the polynomials considered here are always even. It follows that there are no non-trivial symmetric, isotropic tensors of odd order, and that every symmetric isotropic tensor of even order $k$ is a scalar multiple of $\delta^{\odot k}$. This completes the proof.\qed
\medskip
\noindent Given the tensor $T^{i_1...i_{k+2}}$, define the contraction $\delta\llcorner T$ by
$$
(\delta\llcorner T)^{i_1...i_k} = \delta_{pq}T^{i_1...i_kpq}.\eqnum{\nexteqnno[Contraction]}
$$
\proclaim{Lemma \nextprocno}
\noindent For any symmetric tensor $T$ of order $k$,
$$
\delta\llcorner(\delta\odot T)=(m+2k+1)T + \delta\odot(\delta\llcorner T).
\eqnum{\nexteqnno[ContractionOfProductWithDelta]}
$$
\endproclaim
\proclabel{ContractionOfProductWithDelta}
\proof Observe that
$$
(\delta\odot T)^{i_1...i_{k+2}} = \sum_{1\leq p<q\leq k+2}\delta^{i_pi_q}T^{i_1...i_{p-1}i_{p+1}...i_{q-1}i_{q+1}...i_{k+2}}.
$$
Thus
$$\eqalign{
\delta_{i_{k+1}i_{k+2}}(\delta\odot T)^{i_1...i_{k+2}}
&=\delta_{i_{k+1}i_{k+2}}\delta^{i_{k+1}i_{k+2}}T^{i_1...i_k}\cr
&\qquad + \delta_{i_{k+1}i_{k+2}}\sum_{1\leq p\leq k}\delta^{i_pi_{k+2}}T^{i_1...i_{p-1}i_{p+1}...i_{k+1}}\cr
&\qquad + \delta_{i_{k+1}i_{k+2}}\sum_{1\leq p\leq k}\delta^{i_pi_{k+1}}T^{i_1...i_{p-1}i_{p+1}...i_ki_{k+2}}\cr
&\qquad + \sum_{1\leq p,q\leq k}\delta^{i_pi_q}\left(\delta_{i_{k+1}i_{k+2}}T^{i_1...i_{p-1}i_{p+1}...i_{q-1}i_{q+1}...i_{k+2}}\right)\cr
&=\left[(m+1)T + 2kT + \delta\odot(\delta\llcorner T)\right]^{i_1...i_k},}
$$
and the result follows.\qed
\proclaim{Lemma \nextprocno}
\noindent For all $k$,
$$
\delta\llcorner\delta^{\odot k} = k(m + 2k - 1)\delta^{\odot (k-1)}.
\eqnum{\nexteqnno[ContractionOfPowerWithDelta]}
$$
\endproclaim
\proclabel{ContractionOfPowerWithDelta}
\proof We proceed by induction. First observe that $\delta\llcorner\delta=(m+1)$. Next, suppose that it holds for $k$, then, by \eqnref{ContractionOfProductWithDelta} and the inductive hypothesis,
$$\eqalign{
\delta\llcorner\delta^{\odot(k+1)} &= \delta\llcorner(\delta\odot\delta^{\odot k})\cr
&=(m+4k+1)\delta^{\odot k} + \delta\odot(\delta\llcorner\delta^{\odot k})\cr
&=\left((m+4k+1) + k(m + 2k - 1)\right)\delta^{\odot k}\cr
&=(k+1)(m + 2(k+1) - 1)\delta^{\odot k},\cr}
$$
and the result follows.\qed
\proclaim{Lemma \nextprocno}
\noindent For all $k$,
$$\eqalign{
\int_{S^m}x^{i_1}...x^{i_{2k}}\opdVol &= \frac{\opVol(S^m)(m-1)!!}{k!(m+2k-1)!!}\delta^{\odot k}, \&\cr
\int_{S^m}x^{i_1}...x^{i_{2k+1}}\opdVol &= 0.\cr}\eqnum{\nexteqnno[IntegralFormulae]}
$$
\endproclaim
\proclabel{IntegralFormulae}
\proof For all $l$, denote
$$
M_l^{i_1...i_l} := \int_{S^m}x^{i_1}...x^{i_l}\opdVol.
$$
Since $M_l$ is symmetric and isotropic, it follows by Lemma \procref{SymmetricIsotropicTensors} that $M_l$ vanishes when $l$ is odd, and when $l=:2k$ is even $M_l=C_k\delta^{\odot k}$ for some constant $C_k$. It remains to show that
$$
C_k = \frac{\opVol(S^m)(m-1)!!}{k!(m+2k-1)!!}
$$
for all $k$. We prove this by induction on $k$. Indeed, $C_0=\opVol(S^m)$. Now suppose that it holds for $k$. Since $\|x\|^2=1$ over $S^m$, for all $k$,
$$
(\delta\llcorner M_{2(k+1)})^{i_1...i_{2k}} = \delta_{i_{2k+1}i_{2k+2}}\int_{S^m}x^{i_1}...x^{i_{2k+2}}\opdVol = \int_{S^m}x^{i_1}...x^{i_{2k}}\opdVol = M_{2k}^{i_1...i_{2k}},
$$
so that, by \eqnref{ContractionOfPowerWithDelta} and the induction hypothesis,
$$
C_{k+1} = \frac{1}{(k+1)(m+2k+1)}C_k = \frac{\opVol(S^m)(m-1)!!}{(k+1)!(m+2k+1)!!},
$$
and the result follows.\qed
\proclaim{Proposition \nextprocno}
\noindent The functions $(x^i)_{1\leq i\leq m+1}$ constitute an orthogonal basis of $\Cal{H}_1$ with respect to the $L^2$ inner product over $S^m$.
\endproclaim
\proclabel{OrthonormalBasis}
\proof These functions trivially constitute a basis of $\Cal{H}_1$. Furthermore, by Lemma \procref{IntegralFormulae}, for all $1\leq i,j\leq m+1$,
$$
\int_{S^m}x^ix^j\opdVol = \frac{\opVol(S^m)}{(m+1)}\delta^{ij},
$$
and orthogonality follows.\qed
\medskip
\noindent Let $\Pi:L^2(S^m)\rightarrow\Cal{H}_1$ be the orthogonal projection.
\noskipproclaim{Proposition \nextprocno}
$$
\Pi\left(\frac{1}{4}\opRic_{ab;c}x^ax^bx^c - \frac{(m+1)}{2(m+3)}S_{;a}x^a\right) = 0,\eqnum{\nexteqnno[ZerothTermInPhi]}
$$
and, for any fixed vector $V$,
$$
\Pi\left(\frac{1}{4}\opRic_{ab;cd}x^ax^bx^cV^d - \frac{(m+1)}{2(m+3)}S_{;ab}x^aV^b\right) = 0.\eqnum{\nexteqnno[ThirdSummand]}
$$
\endproclaim
\proclabel{ThirdSummand}
\remark Observe that \eqnref{ThirdSummand} corresponds to the third summand in the asymptotic series \eqnref{AsymptoticExpansionOfPhi} of $\Phi$.
\medskip
\proof Indeed, bearing in mind Lemma \procref{IntegralFormulae}, for all $1\leq i\leq m+1$,
$$\multiline{
\int_{S^m}\left(\frac{1}{4}\opRic_{ab;c}x^ax^bx^c - \frac{(m+1)}{2(m+3)}S_{;a}x^a\right)x^i\opdVol\cr
\qquad\qquad = \frac{\opVol(S_m)}{4(m+1)(m+3)}\opRic_{ab;c}(\delta^{ab}\delta^{ci} + \delta^{ac}\delta^{bi} + \delta^{ai}\delta^{bc}) - \frac{\opVol(S_m)}{2(m+3)}S_{;a}\delta^{ai}.\cr}
$$
The first relation now follows by the second Bianchi identity and the second follows upon taking its formal derivative.\qed
\noskipproclaim{Proposition \nextprocno}
$$
\Pi\left(\frac{1}{3}\opRic_{ab}x^ax^b\right) = 0,\eqnum{\nexteqnno[FirstTermInPhi]}
$$
and, for any fixed vector $V$,
$$
\Pi\left(\frac{1}{3}\opRic_{ab;c}x^ax^bV^c\right) = 0,\eqnum{\nexteqnno[FourthSummand]}
$$
\endproclaim
\proclabel{FourthSummand}
\remark Observe that \eqnref{FourthSummand} corresponds to the fourth summand in the asymptotic series \eqnref{AsymptoticExpansionOfPhi} of $\Phi$.
\medskip
\proof The first relation follows directly from Lemma \procref{IntegralFormulae} and the second relation follows upon taking the formal derivative.\qed
\newsubhead{Formal Solutions}[FormalSolutions]
\proclaim{Theorem \nextprocno}
\noindent There exist increasing sequences $(C_k)$ of positive constants and $(n_k)$ of positive integers with the property that, for all $s$, there exist canonical sequences $(Y_{k,s})\in C^{1,\alpha}(\Bbb{R},\Bbb{R}^{m+1})$ and $(f_{k,s})\in\hat{C}^{1,\alpha}_{\opin}(\Bbb{R}\times S^m)$ such that, for all $k$,
$$\eqalign{
f_{k,s} &\in C^{1,\alpha}(\Bbb{R},\hat{\Cal{H}}_{n_k}),\cr
\|f_{k,s}\|_{1,\alpha,\opin,s^4} &\leq C_k,\ \&\cr
\|Y_{k,s}\|_{1,\alpha} &\leq C_k,\cr}
\eqnum{\nexteqnno[SequenceProperties]}
$$
and, for all $N$,
$$
\left\|\Phi\left(s,\sum_{k=0}^{N-1}s^k Y_{k,s},\sum_{k=0}^Ns^kf_{k,s}\right)\right\|_{0,\alpha,\opin}\leq C_k s^{N+1}.
\eqnum{\nexteqnno[SequencesYieldFormalSolution]}
$$
\endproclaim
\proclabel{ExistenceOfApproximateSolutions}
\proof We prove this by induction. First define the projection $\Pi:C^{k,\alpha}_\opin(\Bbb{R}\times S^m)\rightarrow\hat{C}^{k,\alpha}(\Bbb{R},\Cal{H}_1)$ by
$$
\Pi(f)(t,x) := \sum_{i=0}^{m+1}\frac{(m+1)}{\opVol(S^m)}\int_{S^m}f(t,x)x^i\opdVol x^i.
$$
That is, for each $t$, $\Pi(f)(t,\cdot)$ is the $L^2$-orthogonal projection of the function $f(t,\cdot)$ onto $\Cal{H}_1$.  Observe, that, for all $l$, for all $f\in C^{1,\alpha}(\Bbb{R},\hat{\Cal{H}}_l)$ and for all $s$,
$$
\Pi Q_s f = 0,\eqnum{\nexteqnno[QMapsToKernelOfPi]}
$$
so that, by Proposition \procref{FourthSummand}, the terms up to order $k$ in the asymptotic expansion of $\Pi\Phi$ only depend on the asymptotic expansions of $f$ and $Y$ up to order $k-2$. Finally, define $\Pi^\perp:=\opId-\Pi$.
\par
Fix $s>0$, and define $f_{0,s}:=-H_s\Phi_0$. By Proposition \procref{FourthSummand}, $\Phi_0\in C^{0,\alpha}(\Bbb{R},\hat{\Cal{H}}_2)$, and since the restriction of $Q_sH_s$ to this space equals the identity, it follows that with $f_{0,s}$ so defined, the term of order $0$ in the asymptotic expansion \eqnref{AsymptoticExpansionOfPhi} of $\Phi$ vanishes. Furthermore, by Proposition \procref{NormOfRestrictionOfGreensOperator}, there exists $C_0>0$ such that $\|f_{0,s}\|_{1,\alpha,\opin,s^4}\leq C_0$ for all $s$. Finally, by Propositions \procref{ThirdSummand} and \procref{FourthSummand}, the terms of order $0$ and $1$ in the asymptotic expansion of $\Pi\Phi$ both vanish.
\par
Now suppose that we have defined $C_0,...,C_k,n_0,...,n_k,f_{0,s},...,f_{k,s},Y_{0,s},...,Y_{k-1,s}$ such that the terms up to order $k$ and $k+1$ in the asymptotic expansions of $\Phi$ and $\Pi\Phi$ respectively all vanish, for all $s$, and for all $0\leq l\leq k$,
$$\eqalign{
f_{l,s} \in C^{1,\alpha}(\Bbb{R},\hat{\Cal{H}}_{n_l}),\cr
\|f_{l,s}\|_{1,\alpha,\opin,s^4} \leq C_l,\cr}
$$
and for all $0\leq l\leq k-1$,
$$
\|Y_{l,s}\|_{1,\alpha} \leq C_l.
$$
Define
$$
Y_{k,s} := -G\circ\Pi\circ\Phi_{k+2}(f_{0,s},...,f_{k,s},s^4\dot{f}_{0,s},...,s^4\dot{f}_{k-2,s},
Y_{0,s},...,Y_{k-1,s},\dot{Y}_{0,s},...,\dot{Y}_{k-2,s}),
$$
and define $f_{k+1,s}:=-H_s\Pi^\perp\Psi_{k+1,s}$, where
$$\eqalign{
\Psi_{k+1,s} &= \left(\frac{1}{4}\opRic_{pq;rs}x^px^qx^rY^s_{k-1,s} - \frac{(m+1)}{2(m+3)}S_{;pq}x^pY^q_{k-1,s}\right)-\frac{1}{3}\opRic_{pq;r}x^px^qY^r_{k,s}\cr
&\qquad+\Phi_{k+1}(f_{0,s},...,f_{k-1,s},s^4\dot{f}_{0,s},...,s^4\dot{f}_{k-3,s},
Y_{0,s},...,Y_{k-2,2},\dot{Y}_{0,s},...,\dot{Y}_{k-3,s}).\cr}
$$
Since $\Phi_{k+1}$ is a curvature polynomial, and since $f_l$ takes values in $\hat{\Cal{H}}_{n_l}$ for all $0\leq l\leq k$, there exists $n_{k+1}\geq n_k$ such that $\Pi^\perp\Psi_{k+1,s}(t,\cdot)$ is an element of $\hat{\Cal{H}}_{n_{k+1}}$ for all $s$ and for all $t$. By hypothesis, the term of order $k+1$ in the asymptotic expansion of $\Pi\Phi$ vanishes, and so, since the restriction of $Q_sH_s$ to $C^{0,\alpha}(\Bbb{R},\hat{\Cal{H}}_{n_{k+1}})$ equals the identity, with $f_{k+1,s}$ so defined, the term of order $k+1$ in the asymptotic expansion of $\Phi$ vanishes. Finally, observe that the function $\Phi_{k+1}(\cdot,...,\cdot)$ is bounded, and since its derivatives are uniformly bounded in $t$, it is uniformly Lipschitz. There therefore exists $B>0$ such that, for all $s$,
$$
\|\Phi_{k+1,s}\|_{0,\alpha,\opin} \leq B,
$$
and by Proposition \procref{NormOfRestrictionOfGreensOperator}, there exists $C_{k+1}>C_k$ such that, for all $s$,
$$
\|f_{k+1,s}\|_{1,\alpha,\opin,s^4} \leq C_{k+1}.
$$
In like mannner, since $PG$ equals the identity, by Propositions \procref{ThirdSummand} and \procref{FourthSummand}, with $Y_{k,s}$ so defined, the term of order $k+2$ in the asymptotic expansion of $\Pi\Phi$ vanishes. Furthermore, upon increasing $C_k$ if necessary, we may suppose that, for all $s$,
$$
\|Y_{k,s}\|_{1,\alpha}\leq C_k.
$$
\par
We have therefore constructed sequences $(C_k)$, $(n_k)$, $(Y_{k,s})$ and $(f_{k,s})$ satisfying the conclusions of the theorem such that,
$$
\Phi\left(s,\sum_{k=0}^\infty s^k Y_{k,s},\sum_{k=0}^\infty s^k f_{k,s}\right) \sim 0.
$$
Observe that the partial sum of $\Phi$ up to order $N$ only involves terms up to order $N-1$ in $Y$. Furthermore, the time-derivative of $f$ only ever appears together with a factor of $s^4$. Thus, for all $N\geq 0$, there exists a smooth function $R_N$ with uniformly bounded derivatives such that for all $s$ and for all $(t,x)$,
$$\multiline{
\Phi\left(s,\sum_{k=0}^{N-1}s^kY_{k,s,t},\sum_{k=0}^N s^k f_{k,s,t,x}\right) \cr
\qquad\qquad = s^{N+1}R_N(s,Y_{0,s,t},...,Y_{N-1,s,t},\dot{Y}_{0,s,t},...,\dot{Y}_{N-1,s,t},\cr
\qquad\qquad\qquad\qquad f_{N,s,t,x},...,f_{N,s,t,x},s^4\dot{f}_{0,s,t,x},...,s^4\dot{f}_{N,s,t,x}).\cr}
$$
The function $R_N$ is bounded, and since its derivatives are uniformly bounded in $t$, it is uniformly Lipschitz. Thus, upon increasing $C_k$ if necessary, it follows as before that
$$
\left\|\Phi\left(s,\sum_{k=0}^{N-1}s^kY_{k,s,t},\sum_{k=0}^N s^k f_{k,s,t,x}\right)\right\|_{0,\alpha,\opin} \leq C_ks^{N+1},
$$
as desired.\qed
\newsubhead{Exact Solutions}[ExactSolutions]
We recall the classical inverse function theorem (c.f. \cite{Rudin}).
\proclaim{Theorem \nextprocno, {\bf Inverse Function Theorem}}
\noindent Let $E$ and $F$ be Banach spaces. Let $\Omega$ be a neighbourhood of $0$ in $E$. Let $\Phi:\Omega\rightarrow F$ be a $C^2$ mapping. Suppose that there exists $A,B>0$ such that
$$
\|D\Phi(0)^{-1}\| \leq A,\qquad \|D^2\Phi(x)\| \leq B\ \forall\ x\in\Omega.
$$
If $\epsilon:=\|\Phi(0)\|<1/4A^2B$, and if $B_{2A\epsilon}(0)\subseteq\Omega$, then there exists a unique point $x\in B_{2A\epsilon}(0)$ such that $\Phi(x)=0$.
\endproclaim
We now obtain existence.
\proclaim{Theorem \nextprocno}
\noindent For all sufficiently small $s$, there exist canonical functions $Y_s\in C^{1,\alpha}(\Bbb{R},\Bbb{R}^{m+1})$ and $f_s\in C^{1,\alpha}_\opin(\Bbb{R}\times S^n)$ such that $\Phi(s,f_s,Y_s)=0$. Furthermore, there exists a sequence $(C_k)$ of positive numbers such that if $(Y_{k,s})$ and $(f_{k,s})$ are as in Theorem \procref{ExistenceOfApproximateSolutions}, then, for all $N$
$$
\left\|Y_s - \sum_{k=0}^N s^k Y_{k,s}\right\|_{1,\alpha},\left\|f_s - \sum_{k=0}^N s^k f_{k,s}\right\|_{1,\alpha,\opin,s^4} \leq C_Ns^{N+1}.
$$
\endproclaim
\proclabel{ThmExistence}
\proof Let $\Pi$ and $\Pi^\perp$ be as in the proof of Theorem \procref{ExistenceOfApproximateSolutions}. Define the mapping $\Psi:]0,\infty[\times C^{1,\alpha}(\Bbb{R},\Bbb{R}^{m+1})\times\hat{C}^{1,\alpha}_\opin(\Bbb{R}\times S^m)\rightarrow C^{0,\alpha}(\Bbb{R},\Bbb{R}^{m+1})\times\hat{C}^{0,\alpha}_{\opin}(\Bbb{R}\times S^m)$ by
$$
\Psi(s,Y,f) := (s^{-2}\Pi\circ\Phi(s,Y,f),\Pi^\perp\circ\Phi(s,Y,f)).
$$
Consider the asymptotic series \eqnref{AsymptoticExpansionOfPhi} for $\Phi$ up to order $2$ in $s$. Substituting $f_{0,s}=f$, $f_{1,s}=f_{2,s}=0$, $Y_{0,s}=Y$ and $Y_{1,s}=0$, yields
$$\eqalign{
s^{-2}(\Pi\circ\Phi)(s,Y,f) &= PY + (\Pi\circ R_1)(f,s^4\dot{f}) + s(\Pi\circ R_2)(s,f,s^4\dot{f},\dot{Y}),\cr
(\Pi^\perp\circ\Phi)(s,Y,f) &= Q_sf - \frac{1}{3}\opRic_{pq}x^px^q + sR_3(s,f,s^4\dot{f},Y,\dot{Y}),\cr}
$$
for functions $R_1$, $R_2$ and $R_3$ which are smooth at $s=0$. Differentiating with respect to $Y$ and $f$, it follows that
$$
D\Psi(s,Y,f) = \pmatrix P\hfill& A(s,Y,f)\hfill\cr 0\hfill& Q_s\hfill\cr\endpmatrix + sB(s,Y,f),\eqnum{\nexteqnno[SingularDerivative]}
$$
where, for all $R>0$, there exists $\epsilon, C>0$ such that if $s<\epsilon$ and if $\|Y\|_{1,\alpha}+\|f\|_{1,\alpha,\opin}\leq R$, then $\|A(s,Y,f)\|_{0,\alpha,\opin},\|B(s,Y,f)\|_{0,\alpha,\opin}\leq C$. In particular, by Lemma \procref{InfiniteDimensionalFactor}, we may suppose that $D\Psi(s,Y,f)$ is invertible with $\|D\Psi(s,Y,f)\|\leq Cs^{-\alpha}$. Furthermore, we may likewise suppose that for all such $s$, $Y$ and $f$, $D^2\Psi(s,Y,f)\leq C$.
\par
Let $(C_k)$, $(Y_{k,s})$ and $(f_{k,s})$ be as in Theorem \procref{ExistenceOfApproximateSolutions}. Upon reducing $\epsilon$ if necessary, we may suppose that, for all $s<\epsilon$,
$$
\Phi(s,Y_0,f_0 + sf_1) \leq \frac{s^{2\alpha}}{4C^3},
$$
and it follows by the inverse function theorem that for all such $s$, there exists a unique pair $(Y,f)$ such that $\|Y_s\|_{1,\alpha}+\|f_s\|_{1,\alpha,\opin,s^4}<s^\alpha/2C^2$ and $\Phi(s,Y,f)=0$. Now fix $N>0$. Upon reducing $\epsilon$ further if necessary, we may suppose that for all $s<\epsilon$
$$
\Phi\left(s,\sum_{k=0}^{N-1}s^kY_k,\sum_{k=0}^Ns^kf_k\right) \leq Cs^{N+1} < \frac{s^{2\alpha}}{4C^3},
$$
and it follows by the inverse function theorem again there for all such $s$, there exists a unique pair $(Y',f')$ such that $\|Y'_s\|_{1,\alpha} + \|f'_s\|_{1,\alpha,\opin,t^4}< 2C^2s^{N+1-\alpha} < s^\alpha/2C^2$ and $\Phi(s,Y',f')=0$. By uniqueness, $Y'=Y$ and $f'=f$. It follows that for all $N>0$, there exists $\epsilon,C>0$ such that for $s<\epsilon$,
$$
\left\|Y_s - \sum_{k=0}^{N_1}s^kY_k\right\|_{1,\alpha},\left\|f_s - \sum_{k=0}^Ns^k f_k\right\|_{1,\alpha,\opin,s^4} \leq Cs^{N+1-\alpha}.
$$
The result follows.\qed
\goodbreak
\newhead{Bibliography}[Bibliography]
{\leftskip = 5ex \parindent = -5ex
\leavevmode\hbox to 4ex{\hfil \cite{BrendleEichmair}}\hskip 1ex{Brendle S., Eichmair M., Large outlying stable constant mean curvature spheres in initial data sets, {\sl Invent. Math.}, {\bf 197}, (2014), no. 3, 663--682}
\medskip
\leavevmode\hbox to 4ex{\hfil \cite{Chavel}}\hskip 1ex{Chavel I., {\sl Riemannian geometry}, Cambridge Studies in Advanced Mathematics, {\bf 98}, Cambridge University Press, Cambridge, (2006)}
\medskip
\leavevmode\hbox to 4ex{\hfil \cite{EichmairMetzgerI}}\hskip 1ex{Eichmair M., Metzger J., On large volume preserving stable (CMC) surfaces in initial data sets, {\sl J. Diff. Geom}, {\bf 91}, (2012), no. 1, 81--102}
\medskip
\leavevmode\hbox to 4ex{\hfil \cite{EichmairMetzgerII}}\hskip 1ex{Eichmair M., Metzger J., Large isoperimetric surfaces in initial data sets, {\sl J. Diff. Geom.}, {\bf 94}, (2013), no. 1, 159--186}
\medskip
\leavevmode\hbox to 4ex{\hfil \cite{EichmairMetzgerIII}}\hskip 1ex{Eichmair M., Metzger J., Unique isoperimetric foliations of asymptotically flat manifolds in all dimensions, {\sl Invent. Math.}, {\bf 194}, (2013), no. 3, 591--630}
\medskip
\leavevmode\hbox to 4ex{\hfil \cite{Krylov}}\hskip 1ex{Krylov N. V., {\sl Lectures on elliptic and parabolic equations in Sobolev spaces}, Graduate Studies in Mathematics, {\bf 96}, American Mathematical Society, Providence, RI, (2008)}
\medskip
\leavevmode\hbox to 4ex{\hfil \cite{MaxNunSmi}}\hskip 1ex{Maximo D., Nunes I., Smith G., Free boundary minimal annuli in convex three-manifolds, arXiv:1312.5392}
\medskip
\leavevmode\hbox to 4ex{\hfil \cite{NardulliI}}\hskip 1ex{Nardulli S., The isoperimetric profile of a smooth riemannian manifold for small volumes, {\sl Ann. Global Anal. Geom.}, {\bf 36}, (2009), no. 2, 111--131}
\medskip
\leavevmode\hbox to 4ex{\hfil \cite{NardulliII}}\hskip 1ex{Nardulli S., The isoperimetric profile of a non-compact riemannian manifold for small volumes, {\sl Calc. Var. Partial Differential Equations}, {\bf 49}, (2014), no. 1-2, 173--195}
\medskip
\leavevmode\hbox to 4ex{\hfil \cite{Salamon}}\hskip 1ex{Robbin J. W., Salamon D. A., The spectral flow and the Maslov index, {\sl Bull. London Math. Soc.}, {\bf 27}, (1995), 1--33}
\medskip
\leavevmode\hbox to 4ex{\hfil \cite{RosSmi}}\hskip 1ex{Rosenberg H., Smith G., Degree theory of immersed hypersurfaces, arXiv:1010.1879}
\medskip
\leavevmode\hbox to 4ex{\hfil \cite{Rudin}}\hskip 1ex{Rudin W., {\sl Principles of mathematical analysis}, International Series in Pure and Applied Mathematics, McGraw-Hill Book Co., New York, Auckland, D\"usseldorf, (1976)}
\medskip
\leavevmode\hbox to 4ex{\hfil \cite{Schwarz}}\hskip 1ex{Schwarz M., {\sl Morse homology}, Progress in Mathematics, {\bf 111}, Birkh\"auser Verlag, Basel, 1993}
\medskip
\leavevmode\hbox to 4ex{\hfil \cite{SmiCCH}}\hskip 1ex{Smith G., Constant curvature hypersurfaces and the Euler characteristic, arXiv: 1103.3235}
\medskip
\leavevmode\hbox to 4ex{\hfil \cite{SmiEMCFI}}\hskip 1ex{Smith G., Eternal forced mean curvature flows I: a compactness result, {\sl Geom. Dedicata}, {\bf 176}, (2015), 11--29}
\medskip
\leavevmode\hbox to 4ex{\hfil \cite{White}}\hskip 1ex{White B., The space of minimal submanifolds for varying riemannian metrics, {\sl Indiana Math. Journal}, {\bf 40}, (1991), no.1, 161--200}
\medskip
\leavevmode\hbox to 4ex{\hfil \cite{Ye}}\hskip 1ex{Ye R., Foliation by constant mean curvature spheres, {\sl Pacific J. Math.}, {\bf 147}, (1991), no. 2, 381--396}%
\par}
%
%
%
%
\enddocument